\documentclass{article}

\usepackage{tsrml_2022}
\usepackage[utf8]{inputenc} 
\usepackage[T1]{fontenc}    
\usepackage{hyperref}       
\usepackage{url}            
\usepackage{booktabs}       
\usepackage{amsfonts}       
\usepackage{nicefrac}       
\usepackage{microtype}      
\usepackage{xcolor}         

\usepackage{amsfonts}       
\usepackage{amssymb}
\usepackage{amsmath}
\usepackage{amsthm}
\usepackage{mathtools}
\usepackage{natbib}

\usepackage{mdframed}

\title{Rieoptax: Riemannian Optimization in JAX}
\usepackage{graphicx, wrapfig}
\usepackage{float}
\usepackage{graphicx}
\usepackage{minted}
\usepackage{booktabs}
\usepackage{multirow}
\usepackage{graphicx}
\usepackage{subfigure}
\newcommand{\norm}[1]{\left\lVert#1\right\rVert}

\theoremstyle{definition}

\usepackage[utf8]{inputenc} 
\usepackage[T1]{fontenc}    
\usepackage{hyperref}       
\usepackage{url}            
\usepackage{booktabs}       
\usepackage{amsfonts}       
\usepackage{nicefrac}       
\usepackage{microtype}      
\usepackage{xcolor}         

\title{Rieoptax: Riemannian Optimization in JAX with Privacy}

\author{%
  David S.~Hippocampus\thanks{Use footnote for providing further information
    about author (webpage, alternative address)---\emph{not} for acknowledging
    funding agencies.} \\
  Department of Computer Science\\
  Cranberry-Lemon University\\
  Pittsburgh, PA 15213 \\
  \texttt{hippo@cs.cranberry-lemon.edu} \\
}

\begin{document}

\maketitle

\begin{abstract}
   We present Rieoptax, an open source Python library for Riemannian optimization in JAX. We show that many differential geometric primitives, such as Riemannian exponential and logarithm maps, are usually faster in Rieoptax than existing frameworks in Python, both on CPU and GPU. We support various range of basic and advanced stochastic optimization solvers like Riemannian stochastic gradient, stochastic variance reduction, and adaptive gradient methods. A distinguishing feature of the proposed toolbox is that we also support differentially private optimization on Riemannian manifolds. 
\end{abstract}

\section{Introduction}
Riemannian geometry is a generalization of the Euclidean geometry \cite{lee2006riemannian, helgason1979differential} to general Riemannian manifolds. It includes several nonlinear spaces such as the set of positive definite matrices \cite{bhatia2009positive,thanwerdas2021n}, Grassmann manifold of subspaces~\cite{edelman1998geometry,bendokat2020grassmann, absil2009optimization}, Stiefel manifold of orthogonal matrices~\cite{edelman1998geometry,absil2009optimization, chakraborty2019statistics},  kendall shape spaces \cite{kendall1984shape,kendall1989survey, miolane2017template}, hyperbolic spaces \cite{ungar2008gyrovector, ungar2008analytic}, and special Euclidean and orthogonal group \cite{selig2005geometric, gallier2020differential}, to name a few.

Optimization with manifold based constraints has become increasingly popular and has been employed in various applications such as low rank matrix completion~\cite{boumal2011rtrmc}, learning taxonomy embeddings~\cite{nickel2017poincare,nickel2018learning},  neural networks~\cite{huang2017riemannian, huang2017deep, huang2018building,nguyen2019neural,qi2021transductive}, density estimation \cite{hosseini2020alternative}, optimal transport \cite{chewi2020gradient,altschuler2021averaging, shi2021coupling,mishra2021manifold,han2022riemannianblockSPD}, shape analysis~\cite{srivastava2010shape,  huang2016riemannian}, and topological dimension reduction~\cite{kachan2020persistent}, among others. 

In addition, privacy preserving machine learning \cite{dwork2006calibrating,dwork2008differential} has become crucial in real applications, which has been generalized to manifold-constrained problems very recently \cite{reimherr2021differential,utpala2022differentially,han2022differentially}. Nevertheless, such a feature is absent in existing Riemannian optimization libraries \citep{boumal2014manopt,bergmann2022manopt,meghwanshi2018mctorch,kochurov2020geoopt,smirnov2021tensorflow,JMLR:v17:16-177,miolane2020geomstats}. 

In this work, we introduce Rieoptax (\textbf{Rie}mannian \textbf{Opt}imization in J\textbf{ax}), an open source Python library for Riemannian optimization in JAX \cite{frostig2018compiling, jax2018github}. 
The proposed library is mainly driven by the needs of efficient implementation of manifold-valued operations and optimization solvers, readily compatible with GPU and even TPU processors as well as the needs of privacy-supported Riemannian optimization. To the best of our knowledge, Rieoptax is the first library to provide privacy guarantees within the Riemannian optimization framework. 

\subsection{Background on Riemannian optimization, privacy, and JAX}

\paragraph{Riemannian optimization.}
Riemannian optimization \cite{absil2009optimization, boumal2022intromanifolds} considers the following problem
\begin{align}
    \min_{w \in \mathcal{M}} f(w), \label{Eq:OPT}
\end{align}
where $f : \mathcal{M} \rightarrow \mathbb{R}$, and $\mathcal{M}$ denotes a Riemannian manifold. 
Instead of considering \eqref{Eq:OPT} as a constrained problem, Riemannian optimization \cite{absil2009optimization, boumal2022intromanifolds} views it as an unconstrained problem on the manifold space. Riemannian (stochastic) gradient descent \cite{zhang2016first,bonnabel2013stochastic} generalizes the Euclidean gradient descent with intrinsic updates on manifold, i.e., $w_{t+1} = {\rm Exp}_{w_t}(- \eta_t \, {\rm grad} f(w_t))$, where ${\rm grad} f(w_t)$ is the Riemannian (stochastic) gradient, ${\rm Exp}_w(\cdot)$ is the Riemannian exponential map at $w$ and $\eta_t$ is the step size. 
Recent years have witnessed significant advancements for Riemannian optimization where more advanced solvers are generalized from the Euclidean space to Riemannian manifolds. These include
variance reduction methods \cite{zhang2016riemannian, sato2019riemannian, kasai2018riemannian, zhou2019faster, hanmomentum2021,han2021improved}, adaptive gradient methods \cite{becigneul2018riemannian, kasai2019riemannian}, accelerated gradient methods \cite{han2022riemannian, liu2017accelerated,ahn2020nesterov,zhang2018estimate, alimisis2020continuous}, quasi-Newton methods \citep{huang2015broyden,qi2010riemannian}, zeroth-order methods \cite{li2022stochastic} and second order methods, such as trust region methods \cite{absil2007trust} and cubic regularized Newton's methods \cite{agarwal2021adaptive}.

\paragraph{Differential privacy on Riemannian manifolds.}  Differential privacy ({DP}) provides a rigorous treatment for data privacy by precisely quantifying the deviation in the model's output distribution under modification of a small number of data points \citep{dwork2006calibrating,dwork2006our,dwork2008differential,dwork2014algorithmic}. Provable guarantees of {DP} coupled with properties like immunity to arbitrary post-processing and graceful composability have made it a de-facto standard of privacy with steadfast adoption in the real applications \citep{erlingsson2014rappor, apple2017learning,ding2017collecting, near2018differential, abowd2018us}. Further, it has been shown empirically that {DP} models resist various kinds of leakage attacks that can cause privacy violations \citep{rahman2018membership, 10.5555/3361338.3361358,sablayrolles2019white,zhu2019deep,balle2022reconstructing}.

Recently, there is a surge of interest on differential privacy over Riemannian manifolds, which has been explored in the context of Fr\'echet mean computation~\citep{reimherr2021differential,utpala2022differentially} and, more generally, empirical risk minimization problems where the parameters are constrained to lie on a Riemannian manifold~\citep{han2022differentially}.

\paragraph{JAX and its ecosystem.} JAX \cite{frostig2018compiling, jax2018github} is recently introduced machine learning framework which support automatic differentiation capabilities \cite{baydin2018automatic} via  \texttt{grad()}. Further some of the distinguishing features of JAX are just-in-time (JIT) compilation using the accelerated linear algebra (XLA) compiler \cite{50530} via \texttt{jit()}, automatic vectorization (batch-level parallelism) support with \texttt{vmap()}, and strong support for parallel computation via \texttt{pmap()}.  All the above transformations can be composed arbitrarily because JAX follows the functional programming paradigm and implements these as pure functions.


Given that JAX has many interesting features, its ecosystem has been constantly expanding in the last couple of years. Examples include neural network modules (Flax \cite{flax2020github}, Haiku \cite{haiku2020github}, Equinox \cite{kidger2021equinox}, Jraph \cite{jraph2020github}, Equivariant-MLP \cite{finzi2021practical}), reinforcement learning agents (Rlax \cite{deepmind2020jax}), Euclidean optimization algorithms (Optax \cite{deepmind2020jax}), federated learning (Fedjax \cite{fedjax2021}), optimal transport toolboxes (Ott \cite{cuturi2022optimal}), sampling algorithms (Blackjax \cite{blackjax2020github}), differential equation solvers (Diffrax \cite{kidger2021on}), rigid body simulators (Brax \cite{brax2021github}), and differentiable physics (Jax-md \cite{jaxmd2020}), among others. 


\subsection{Rieoptax} 

We believe that the proposed framework for Riemannian optimization in JAX is a timely contribution that brings several benefits of JAX and new features (such as privacy support) to the manifold optimization community discussed below.

\begin{itemize}
    \item \textbf{Automatic and efficient vectorization with} \texttt{vmap().} Functions that are written for inputs of size $1$ can be converted to functions that take batch of inputs by wrapping it with \texttt{vmap()}. For example, the function \texttt{def dist(point\_a, point\_b)} for computing distance between a single \texttt{point\_a} and a single \texttt{point\_b} can be converted to function that computes distance between a batch of \texttt{point\_a} and/or a batch \texttt{point\_b} by wrapping \texttt{dist} with \texttt{vmap()} without modifying the \texttt{dist()} function. This is useful in many cases, e.g., Fr\'echet mean computation $\min_{w \in \mathcal{M}} \left\{ \frac{1}{n} \sum_{i=1}^{n} f_{i}(w) := \frac{1}{n} \sum_{i=1}^{n} \text{dist}^2(w,z_i)\right\}$. Furthermore, vectorization with \texttt{vmap()} is usually faster or on par with manual vectorization \cite{jax2018github}. 
    \item \textbf{Per-example gradient clipping.} A key process in differentially private optimization is per-example gradient clipping $ \frac{1}{n}\sum_{i=1}^{n} \text{clip}_{\tau}(\text{grad} f_{i}(w))$ , where $\text{clip}_{\tau}$ ensures norm is atmost $\tau$.  
    Here, the order of operations is important: the gradients are first clipped and then averaged. Popular libraries including Autograd~\cite{maclaurin2015autograd}, Pytorch~\cite{paszke2019pytorch} and Tensorflow~\cite{abadi2016tensorflow} are heavily optimized to directly compute the mean gradient $\frac{1}{n}\sum_{i=1}^{n} \text{grad} f_{i}(w)$ and hence do not expose per-example gradients i.e., $\text{grad} f_{i}(w).$ Hence, one has to resort to ad-hoc techniques \cite{goodfellow2015efficient,  rochette2019efficient, lee2021scaling} or come up with algorithmic modifications \cite{bu2021fast} which inherently have speed versus performance trade-off. JAX, however, offers native support for handling such scenarios and JAX-based differentially private Euclidean optimization methods have been shown to be much faster than their non-JAX counterparts~\citep{subramani2021enabling}. We observe that JAX offer similar benefits for differentially private Riemannian optimization as well. 
    
     \item \textbf{Single Source Multiple Devices (SSMD) paradigm.} JAX follows the SSMD paradigm, and hence, the code written for CPUs can be run on GPU/TPUs without any additional modification. 
     
\end{itemize}

Rieoptax is on \url{https://anonymous.4open.science/r/Rieoptax} for review and will be made public.


\section{Design and Implementation overview}

The package currently implements several commonly used geometries, optimization algorithms and differentially private mechanisms on manifolds. More geometries and advanced solvers will be added in the future.

\subsection{Core}

\begin{itemize}
    \item \texttt{rieoptax.core.ManifoldArray :} Lightweight wrapper of \texttt{jax} device array with \texttt{manifold} attribute and used to model array constrained to manifold
    
    \item \texttt{rieoptax.core.rgrad :} Riemannian gradient operator (which is higher order function like $\texttt{grad}$)
\end{itemize}

\subsection{Geometries}

Geometry module contains manifolds equipped with Riemannian metrics. Each Geometry contains Riemannian inner product \texttt{inp()}, and induced norm \texttt{norm()}, Riemannian exponential \texttt{exp()} and logarithm maps \texttt{log()}, induced Riemannian distance \texttt{dist()}, parallel transport \texttt{pt()}, and transformation from Euclidean gradient to Riemannian gradient \texttt{egrad\_to\_rgrad()}.

Manifolds include symmetric positive definite (SPD) matrices ${\rm SPD}(m) := \{ \mathbf X \in \mathbb R^{m \times m} : \mathbf X = \mathbf X^\top, \mathbf X \succ 0 \}$, hyperbolic space, Grassmann manifold $\mathcal{G}(m,r) := \{ [\mathbf X] : \mathbf X \in \mathbb R^{m \times r}, \mathbf X^\top \mathbf X = \mathbf I \}$ where $[\mathbf X] := \{ \mathbf {XO} : O \in O(r) \}$, $O(r)$ denotes the orthogonal group and hypersphere $\mathcal{S}(d) := \{ \mathbf x \in \mathbb R^d : \mathbf x^\top \mathbf x = 1 \}$. We use $T_x \mathcal M$ to represent the tangent space at $x$ and $\langle u,v \rangle_x$ to represent the Riemannian inner product. For more detailed treatment on these geometries, we refer to \cite{absil2009optimization,boumal2022intromanifolds,ungar2008gyrovector}.


\begin{itemize}
    \item \texttt{rieoptax.geometry.spd.SPDAffineInvariant:} SPD matrices with the affine-invariant metric \cite{pennec2006riemannian}: SPD$(m)$ with $\langle \mathbf U, \mathbf V \rangle_{\mathbf X} = {\rm tr}( \mathbf X^{-1} \mathbf U \mathbf X^{-1} \mathbf V)$ for $\mathbf U, \mathbf V \in T_{\mathbf X} {\rm SPD}(m)$.
    
    \item \texttt{rieoptax.geometry.spd.SPDLogEuclidean:} SPD matrices with the log-Euclidean metric \cite{arsigny2007geometric}: SPD$(m)$ with $\langle \mathbf U, \mathbf V \rangle_{\mathbf X} = {\rm tr} \big({\rm D}_{\mathbf U} {\rm logm}(\mathbf X) {\rm D}_{\mathbf V} {\rm logm}(\mathbf X) \big)$ where ${\rm D}_{\mathbf U} {\rm logm}(\mathbf X)$ is the directional derivative of matrix logarithm at $\mathbf X$ along $\mathbf U$.

    \item \texttt{rieoptax.geometry.hyperbolic.PoincareBall:} Poincare-ball model of Hyperbolic space with Poincare metric \cite{ungar2008gyrovector}, i.e., $\mathbb D(d) := \{ \mathbf x \in \mathbb R^{d} : \mathbf x^\top \mathbf x < 1 \}$ with $\langle \mathbf u, \mathbf v \rangle_{\mathbf x} = 4\mathbf u^\top \mathbf v/(1 - \mathbf x^\top \mathbf x)^2$ for $\mathbf u, \mathbf v \in T_{\mathbf x} \mathbb D(d)$.
    
    \item \texttt{rieoptax.geometry.hyperbolic.LorentzHyperboloid:} Lorentz Hyperboloid model of Hyperbolic space \cite{ungar2008gyrovector}, i.e., $\mathbb{H}(d) = \{ \mathbf{x} \in \mathbb{R}^{d} : \langle \mathbf x, \mathbf x \rangle_{\mathcal{L}}  = -1 \}$ with $\langle \mathbf u, \mathbf v \rangle_{\mathbf x} = \langle \mathbf u, \mathbf v \rangle_{\mathcal{L}}$ for $\mathbf u, \mathbf v \in T_{\mathbf x} \mathbb{H}(d)$, where $\langle \mathbf u, \mathbf v \rangle_{\mathcal{L}} := -u_0 v_0 + u_1 v_1 + \cdots u_{d-1} v_{d-1}$.
    
    \item \texttt{rieoptax.geometry.grassmann.GrassmannCanonicalMetric:} Grassmann manifold with the canonical metric \cite{edelman1998geometry}, i.e., $\mathcal G(m,r)$ with $\langle  \mathbf{U}, \mathbf{V}\rangle_{\mathbf {X}}   = {\rm tr}\big( \mathbf{U}^{T} \mathbf{V} \big)$ for $\mathbf U, \mathbf V \in T_{\mathbf X} \mathcal G(m,r)$.
    
    \item  \texttt{rieoptax.geometry.hypersphere.HypersphereCanonicalMetric:}   Hypersphere manifold which canonical metric \cite{absil2009optimization, boumal2022intromanifolds}, i.e., $\mathcal{S}(d)$ with $\langle \mathbf u, \mathbf v \rangle_{\textbf{x}} = \mathbf u^\top \mathbf v$ for $\mathbf u, \mathbf v \in T_{\mathbf x} \mathcal{S}(d)$.
\end{itemize}

\subsection{Optimizers}
Optimizers module contains Riemannian optimization algorithms. Design of optimizers follows Optax \cite{deepmind2020jax}, which implements every optimizer by chaining of few common transformations. Where every optimizer 
\begin{itemize}
    \item \texttt{riepotax.optimizers.first\_order.rsgd:} Riemannian stochastic gradient descent \cite{bonnabel2013stochastic}.

    \item \texttt{riepotax.optimizers.first\_order.rsvrg:}  Riemannian stochastic variance reduced gradient descent \cite{zhang2016riemannian}.
    
     \item \texttt{riepotax.optimizers.first\_order.rsrg:} Riemannian stochastic recursive gradient descent \cite{kasai2018riemannian}.
    
    \item \texttt{riepotax.optimizers.first\_order.rasa:} Riemannian adaptive stochastic gradient algorithm  \cite{kasai2019riemannian}.
    
    \item \texttt{riepotax.optimizers.zeroth\_order.zo\_rgd:} Zeroth-order Riemannian gradient descent \cite{li2022stochastic}.
\end{itemize}

\subsection{Privancy mechanism}

Mechanism module contains differential private mechanisms on Riemannian manifolds.  

\begin{itemize}
    \item \texttt{rieoptax.mechanism.output\_perturbation.RieLaplaceMechanism:} the Riemannian Laplace mechanism \cite{reimherr2021differential} which is used for privatizing Fr\'echet mean computation. 
    
     \item \texttt{rieoptax.mechanism.output\_perturbation.LogEuclideanMechanism:} the Log-Euclidean mechanism \cite{utpala2022differentially} which is used for differentially private Fr\'echet mean on SPD matrices with log-Euclidean metric.
     
    \item \texttt{rieoptax.mechanism.gradient\_perturbation.DPRGDMechanism:} noise calibration for differentially private Riemannian gradient descent \cite{han2022differentially} based on moments accountant \cite{abadi2016deep} in the \texttt{autodp} library \cite{wang2019subsampled}. 
    
    \item \texttt{rieoptax.mechanism.gradient\_perturbation.DPRSGDMechanism:}  noise calibration for Differentially private Riemannian stochastic gradient descent \cite{han2022differentially} based on moments accountant \cite{abadi2016deep} in \texttt{autodp} library \cite{wang2019subsampled}. 
\end{itemize}



\section{Benchmarking Rieoptax}
\label{sect:benchmark}
In this section, we benchmark the proposed Rieoptax against existing Riemannian optimization libraries in Python. These include Pytorch~\cite{paszke2019pytorch} based Mctorch~\cite{meghwanshi2018mctorch} and Geoopt~\cite{kochurov2020geoopt}, Tensorflow~\cite{smirnov2021tensorflow} based Tensorflow-Riemopt (Tf-Riemopt)~\cite{smirnov2021tensorflow}, Numpy~\cite{harris2020array} based Pymanopt \cite{JMLR:v17:16-177}, and Tensorflow based Geomstats \cite{miolane2020geomstats}. While Geomstats supports Numpy, Pytorch, and Tensorflow as backend, currently only the Tensorflow backend provides support for GPUs. 
Other non-Python based libraries include Manopt \cite{boumal2014manopt} in Matlab and Manopt.jl \cite{bergmann2022manopt} in Julia \cite{bezanson2017julia}.

We benchmark the Riemannian exponential (Exp) and logarithm (Log) maps with the proposed Rieoptax against the aforementioned Python libraries whenever available with \texttt{64bitfloat} precision. For CPU benchmarking, we use the AMD EPYC 7B1 processor with 2 cores and 16{GB} RAM. For GPU benchmarking, we use {CUDA} version $11.4$ on 16{GB} Tesla {P100}.

\begin{figure}[t!]
\centering     
\subfigure[Hypersphere Exp]{\label{fig:HypExp1}\includegraphics[width=30mm]{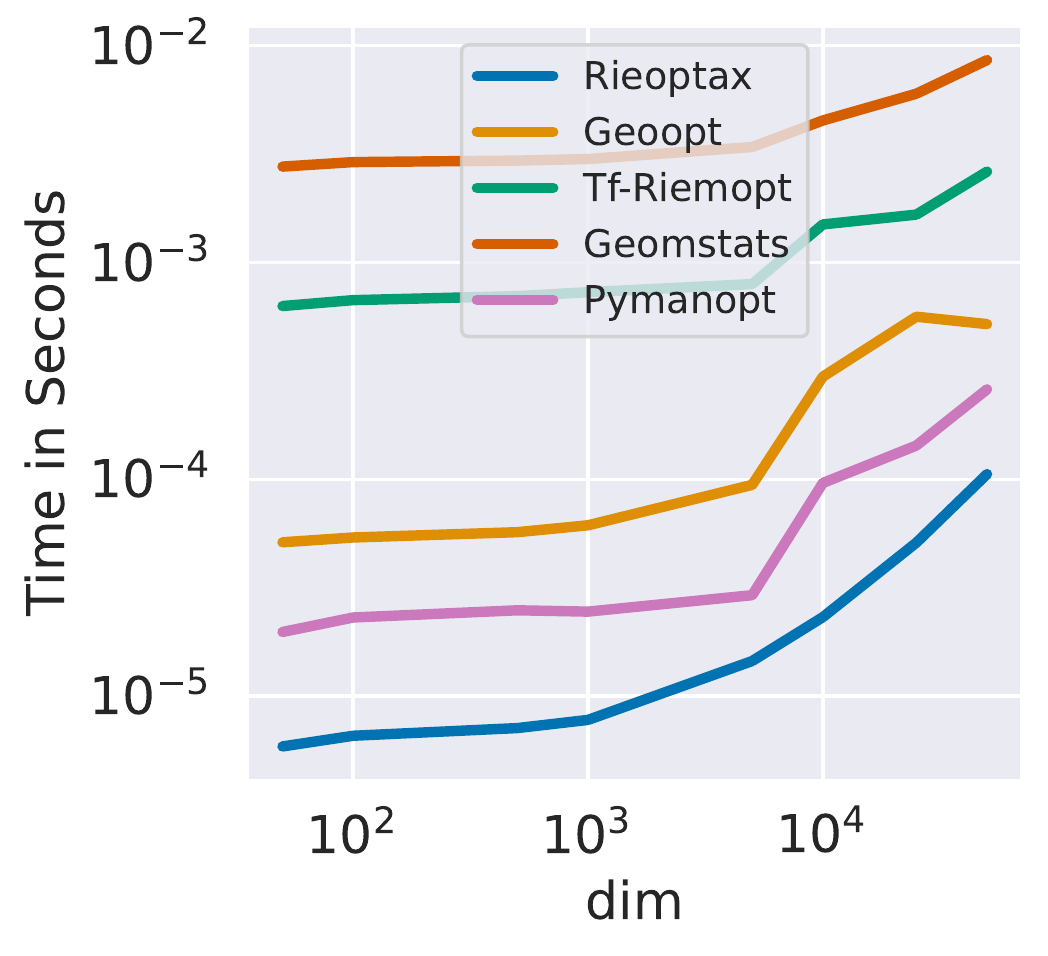}}
\subfigure[Lorentz Exp]{\label{fig:LorExp1}\includegraphics[width=30mm]{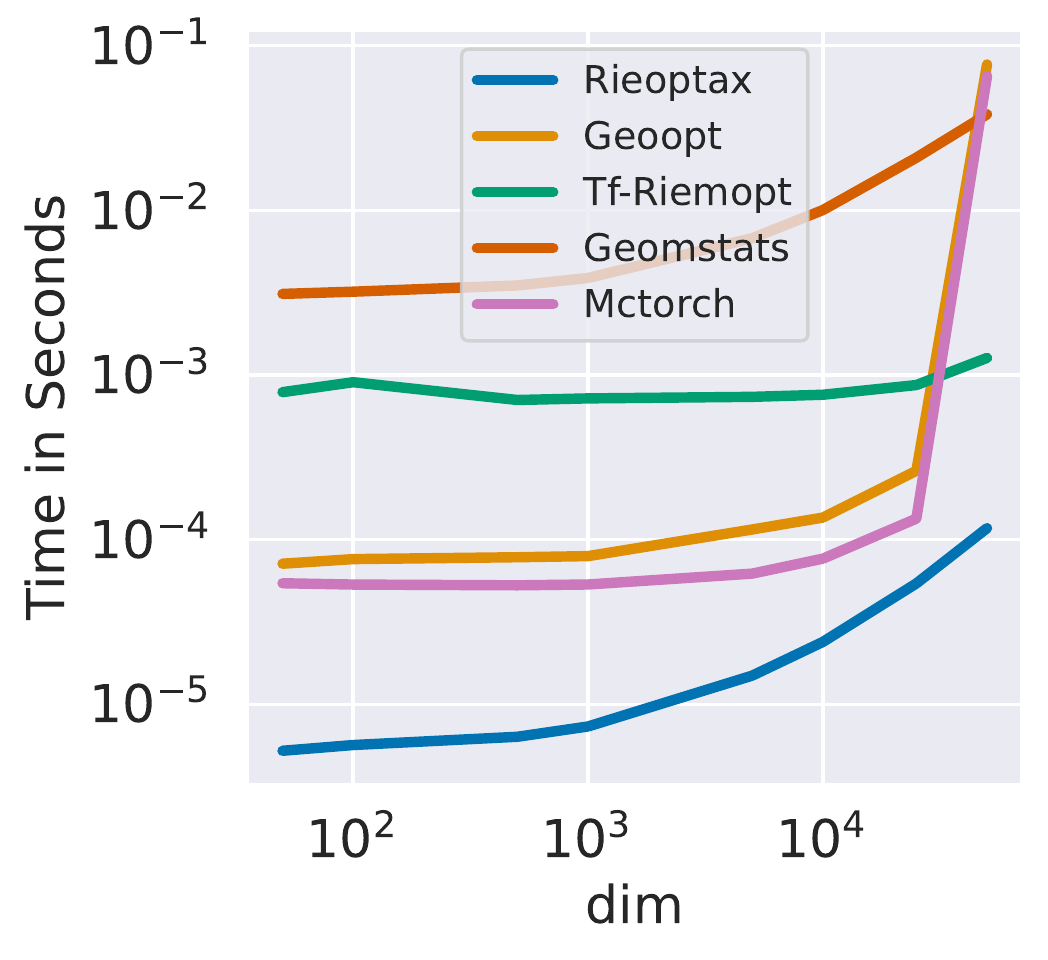}}
\subfigure[Grassmann Exp]{\label{fig:GrassExp1}\includegraphics[width=30mm]{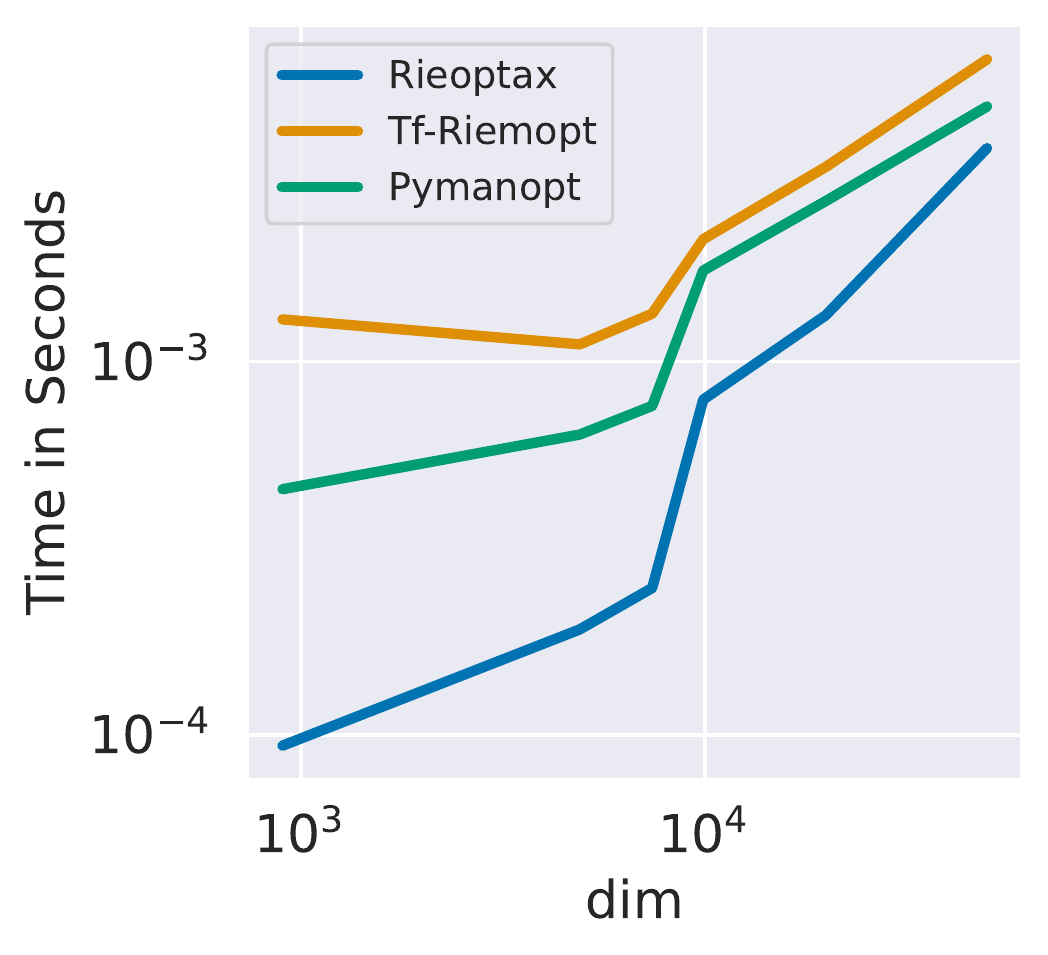}}
\subfigure[SPD Exp]{\label{fig:LorLog1}\includegraphics[width=30mm]{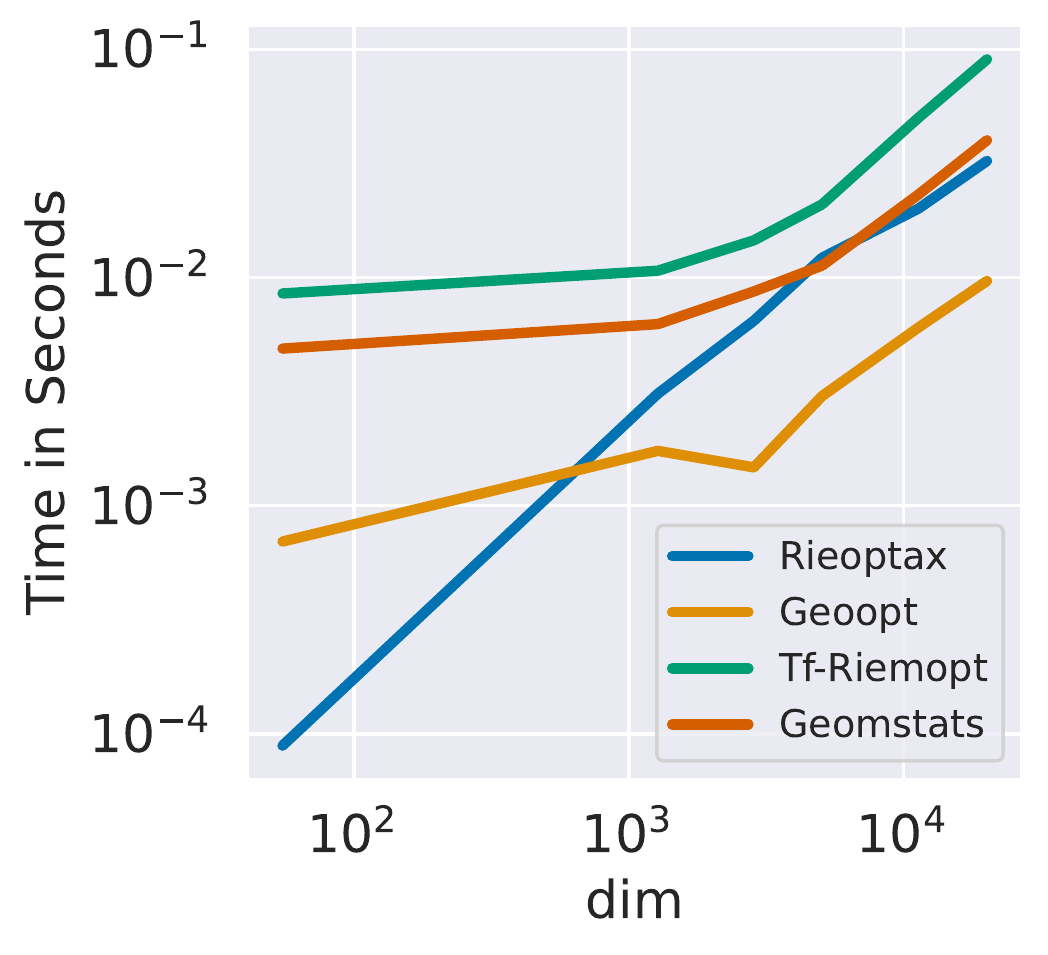}}

\subfigure[Hypersphere Log]{\label{fig:HypLog2}\includegraphics[width=30mm]{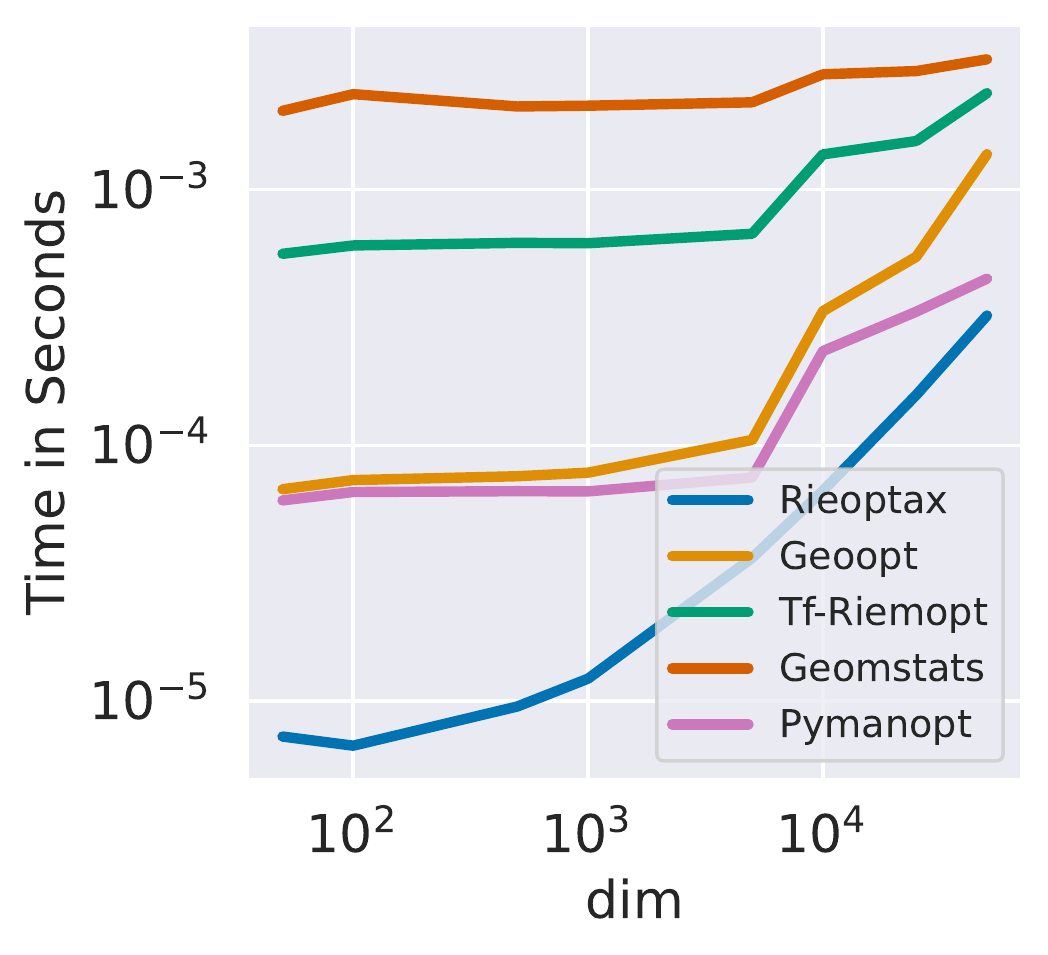}}
\subfigure[Lorentz Log]{\label{fig:LorExp2}\includegraphics[width=30mm]{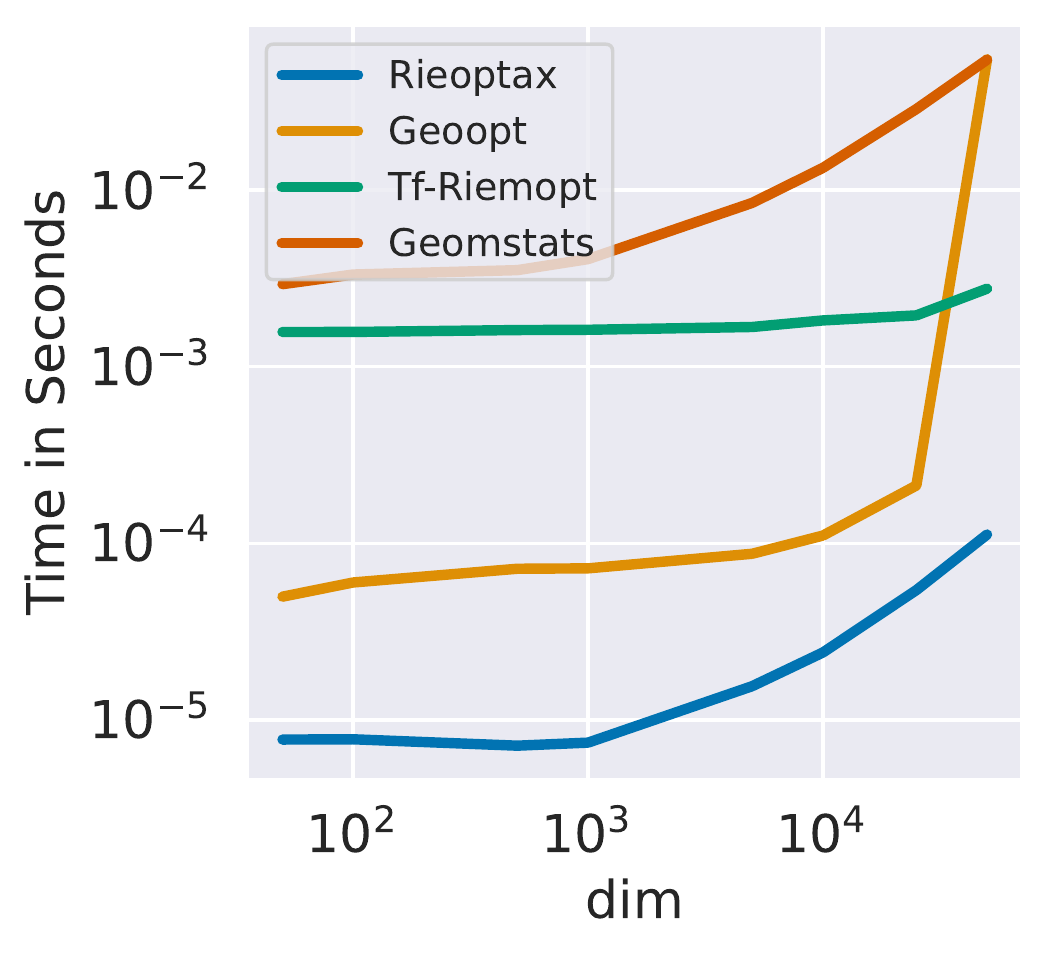}}
\subfigure[Grassmann Log]{\label{fig:GrassLog2}\includegraphics[width=30mm]{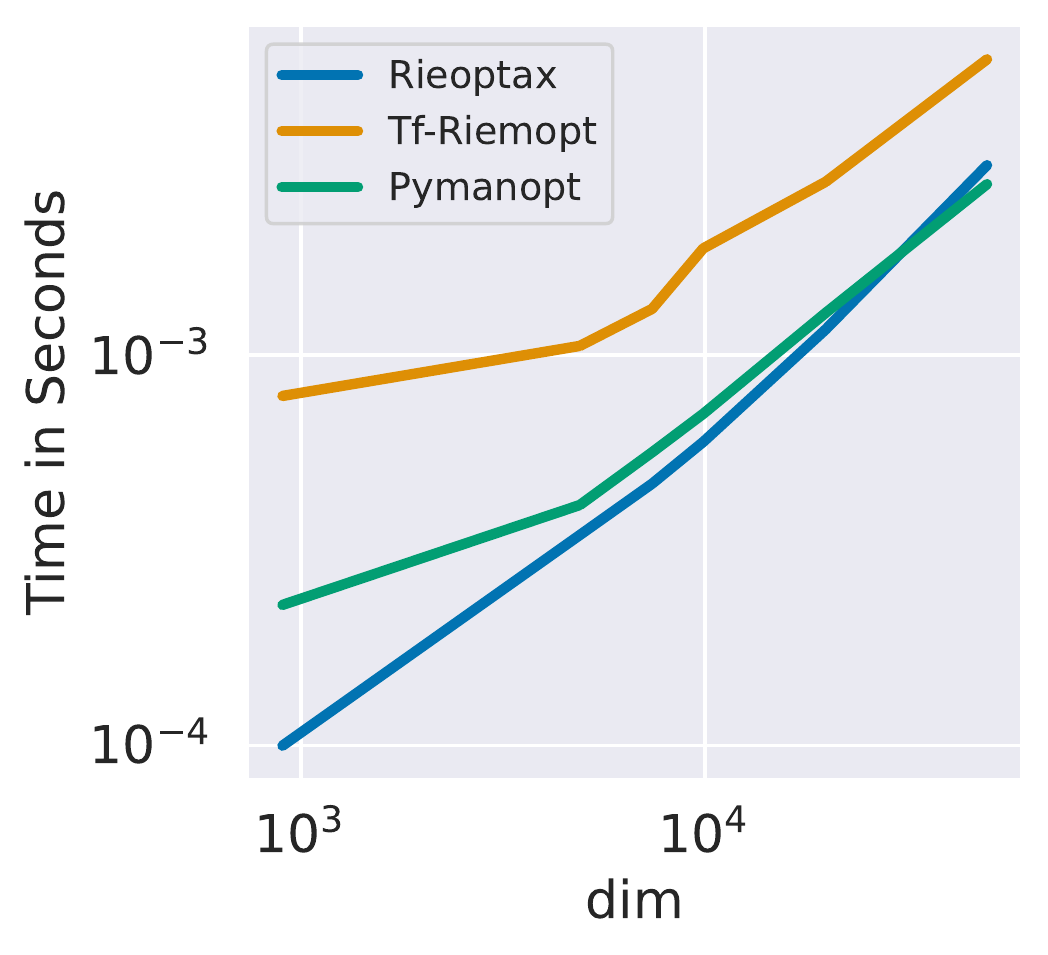}}
\subfigure[SPD  Log]{\label{fig:LorLog2}\includegraphics[width=30mm]{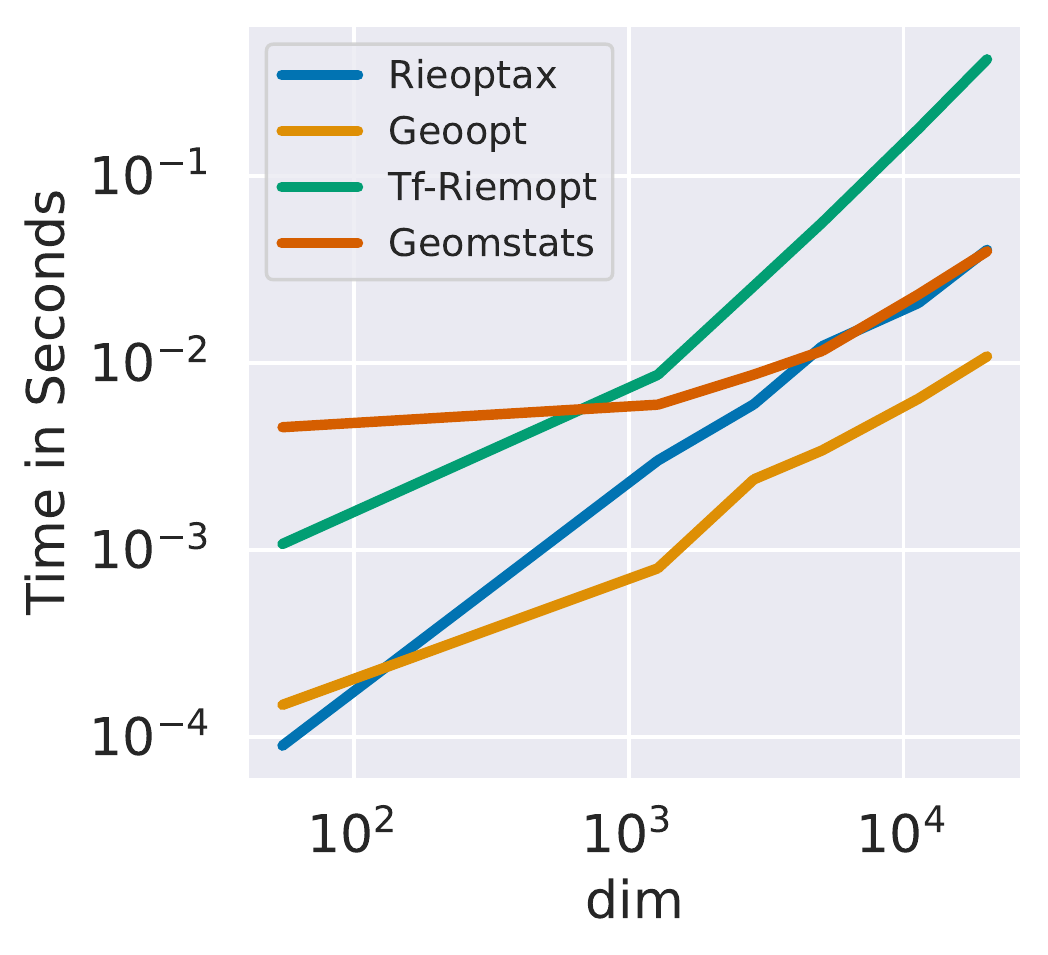}}

\caption{Benchmarking of Geometric Primitives on CPU.}\label{fig:CPUresults}
\end{figure}
\begin{figure}[!]
\centering     
\subfigure[Hypersphere Exp]{\label{fig:HypExp3}\includegraphics[width=30mm]{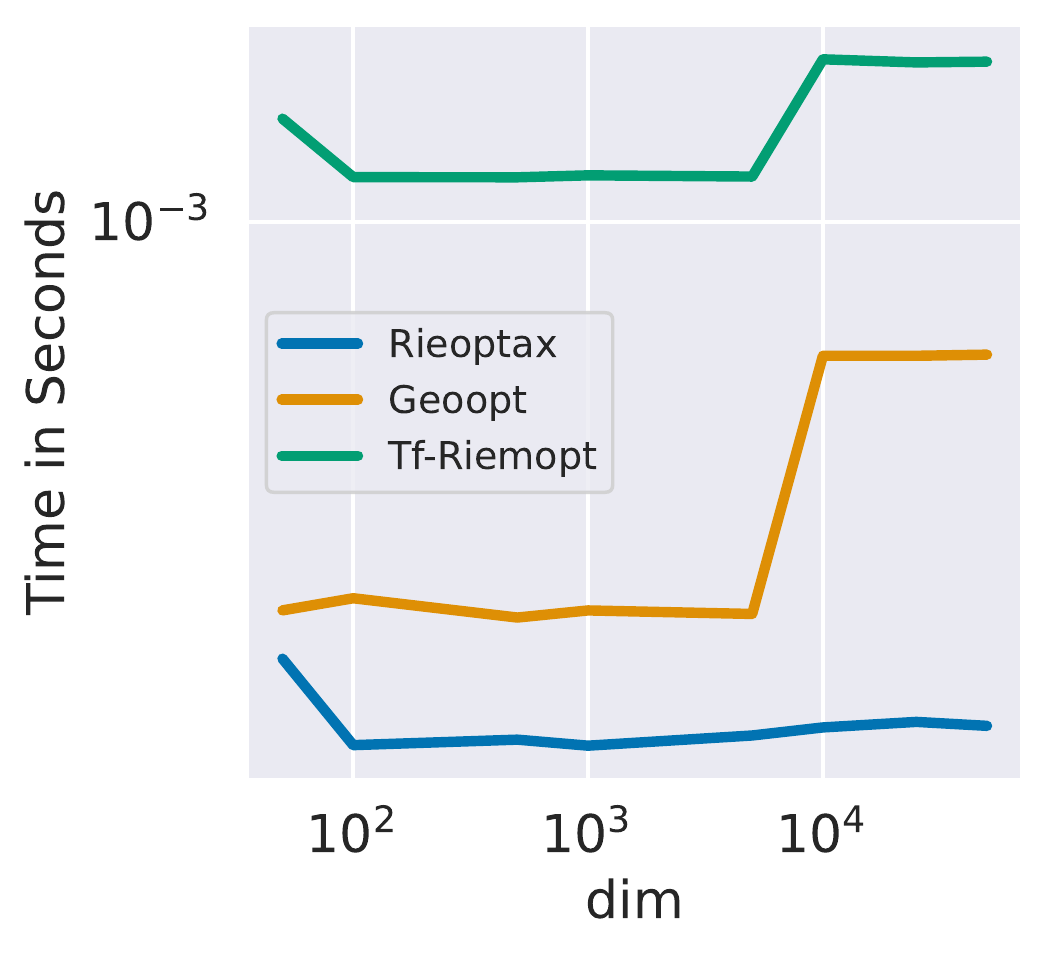}}
\subfigure[Lorentz Exp]{\label{fig:LorExp3}\includegraphics[width=30mm]{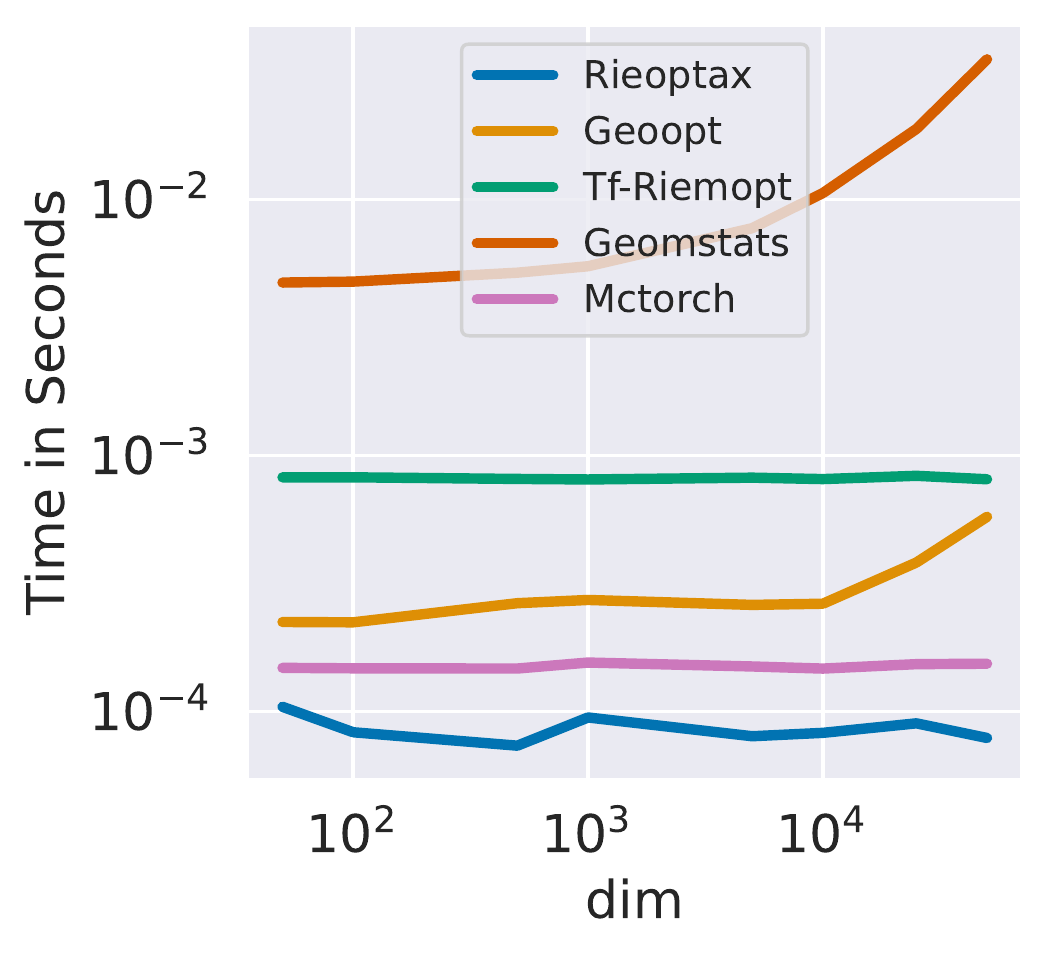}}
\subfigure[Grassmann Exp]{\label{fig:GrassExp3}\includegraphics[width=30mm]{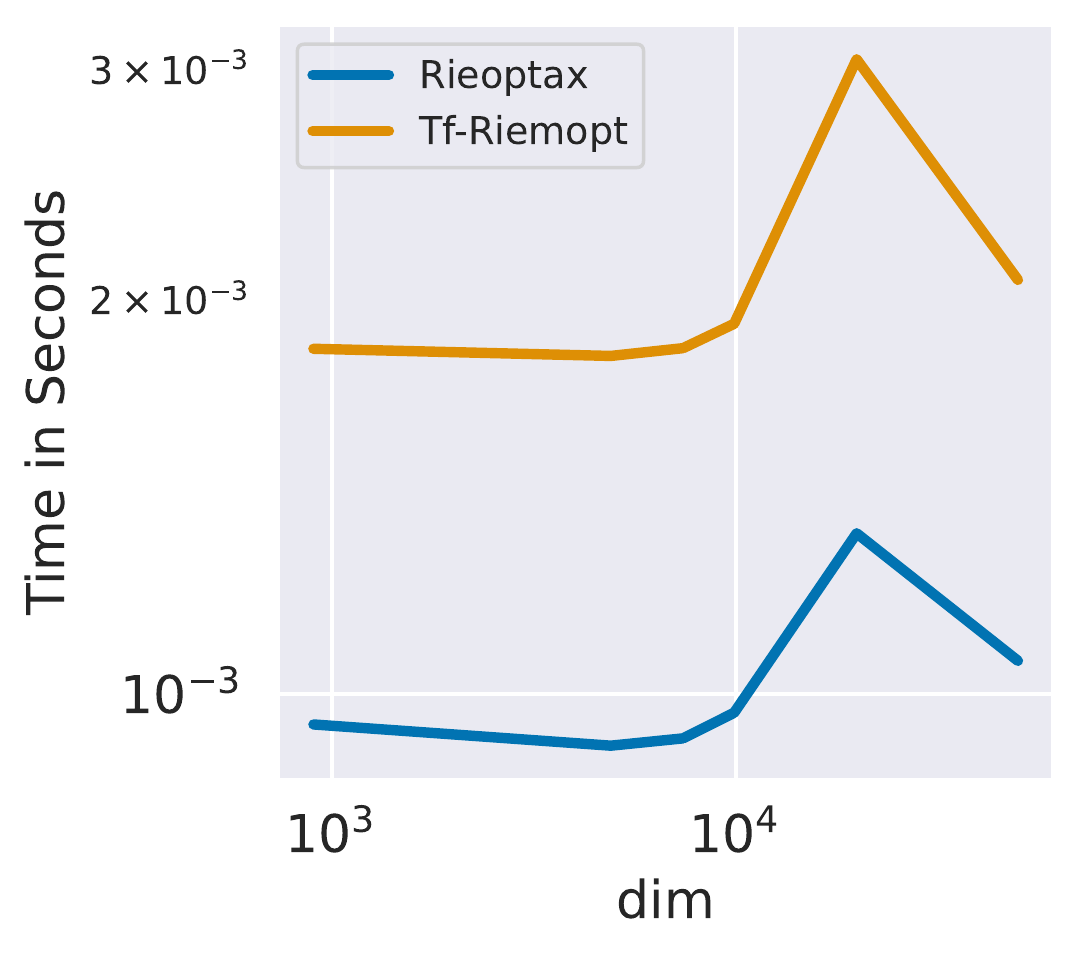}}
\subfigure[SPD Exp]{\label{fig:LorLog3}\includegraphics[width=30mm]{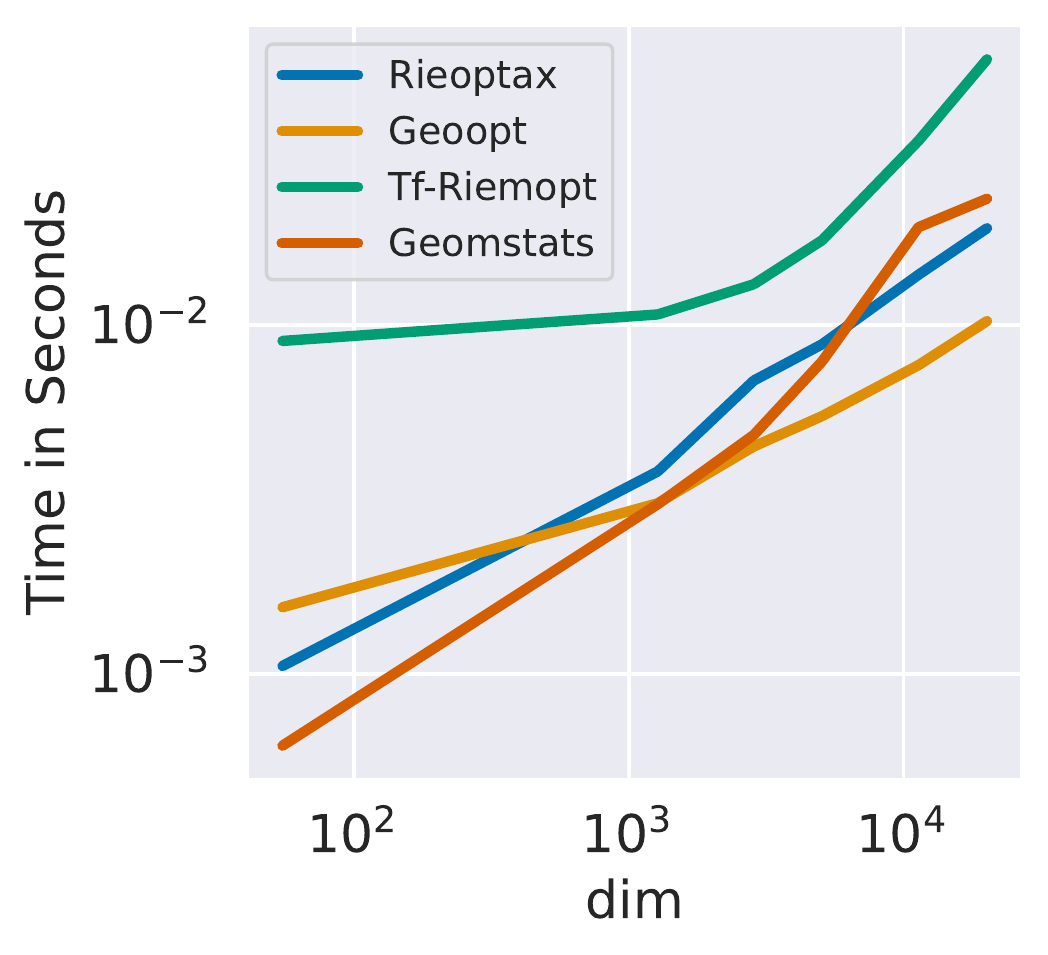}}

\subfigure[Hypersphere Log]{\label{fig:HypExp}\includegraphics[width=30mm]{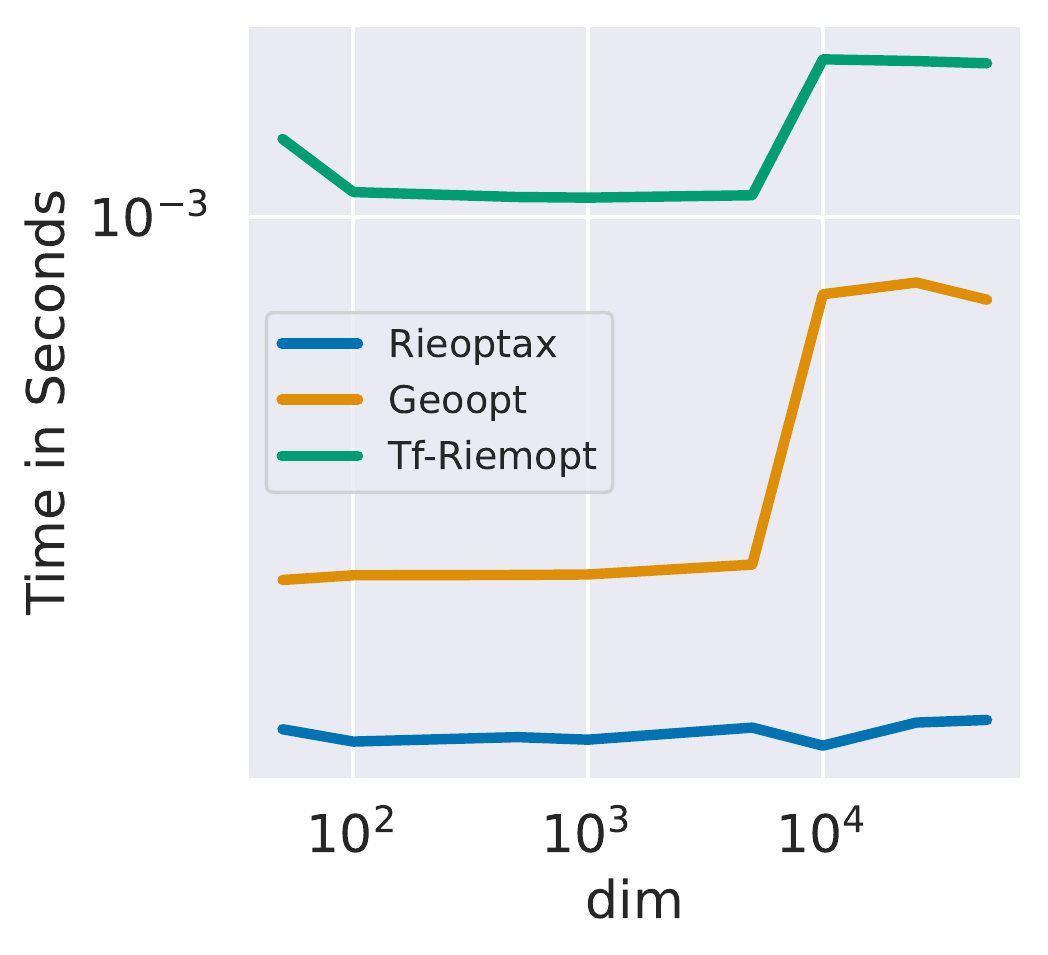}}
\subfigure[Lorentz Log]{\label{fig:LorExp}\includegraphics[width=30mm]{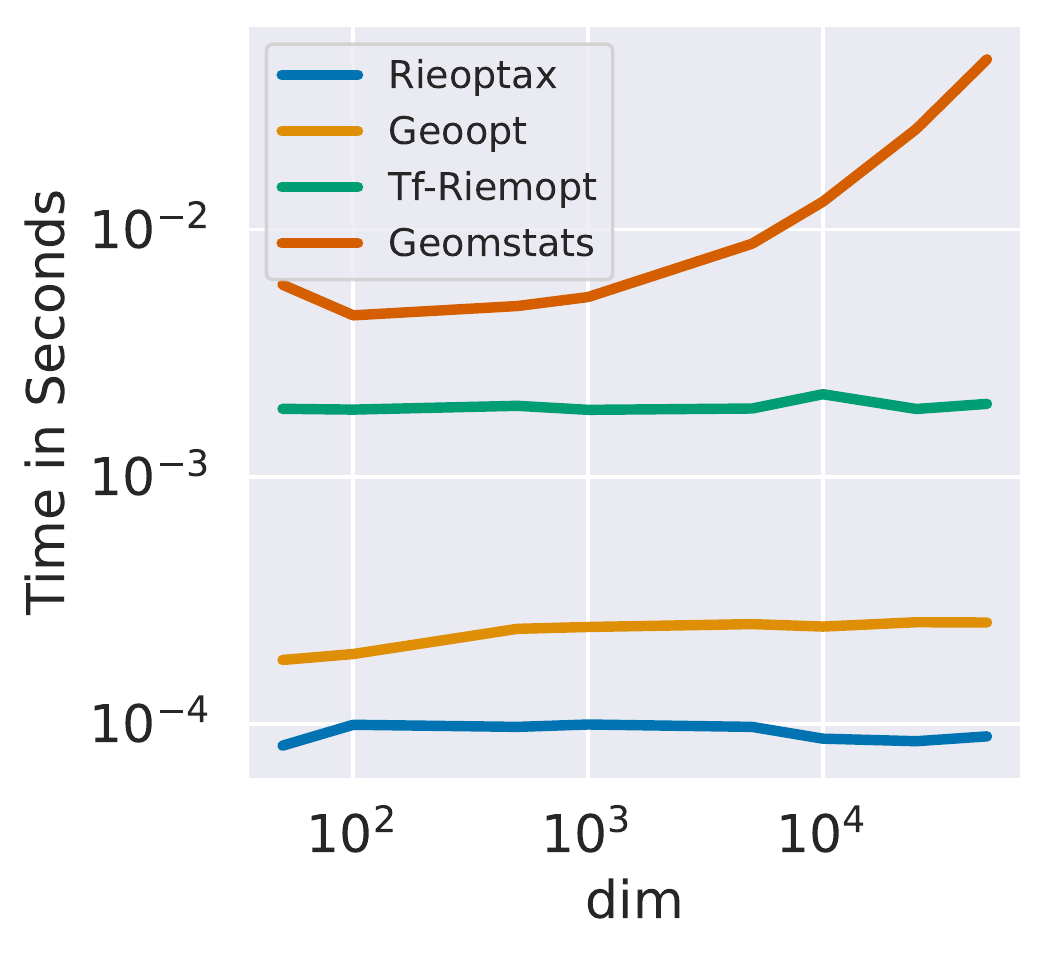}}
\subfigure[Grassmann Log]{\label{fig:GrassExp}\includegraphics[width=30mm]{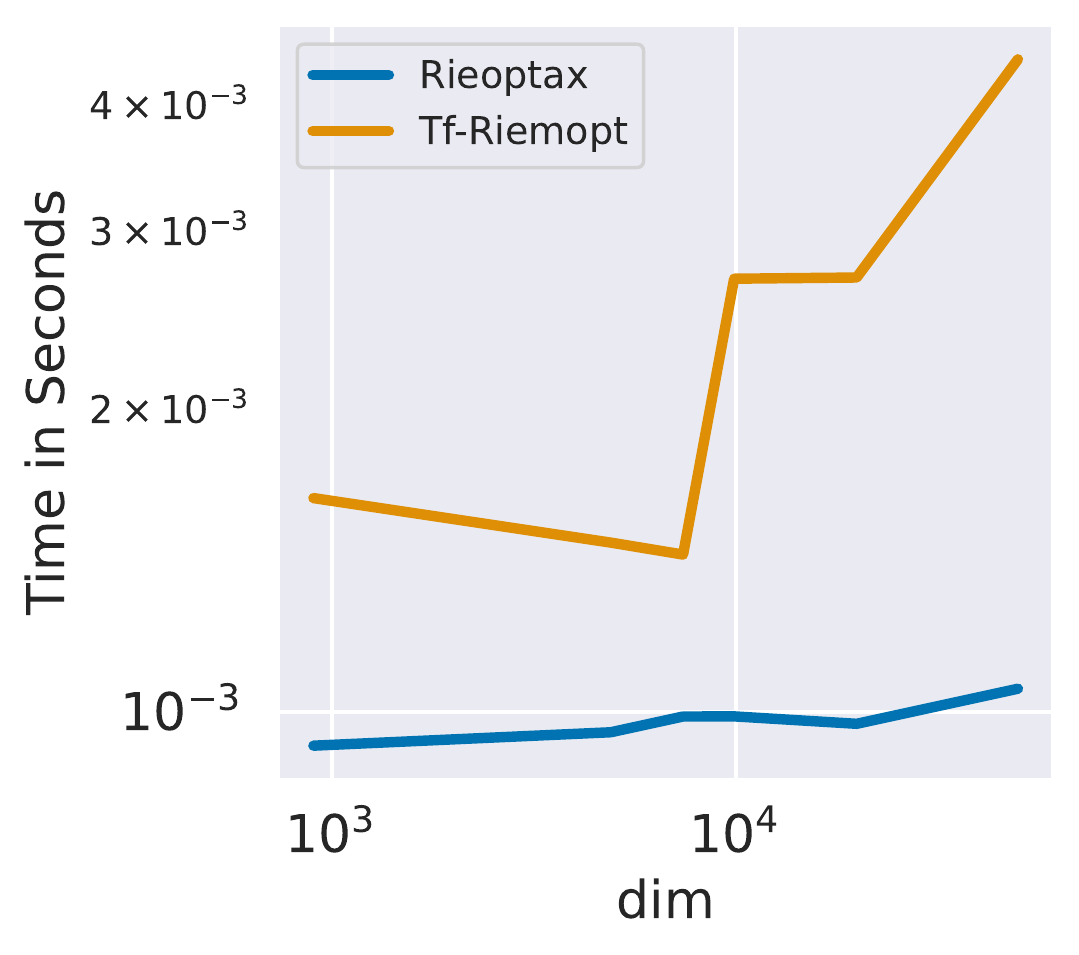}}
\subfigure[SPD Log]{\label{fig:LorLog}\includegraphics[width=30mm]{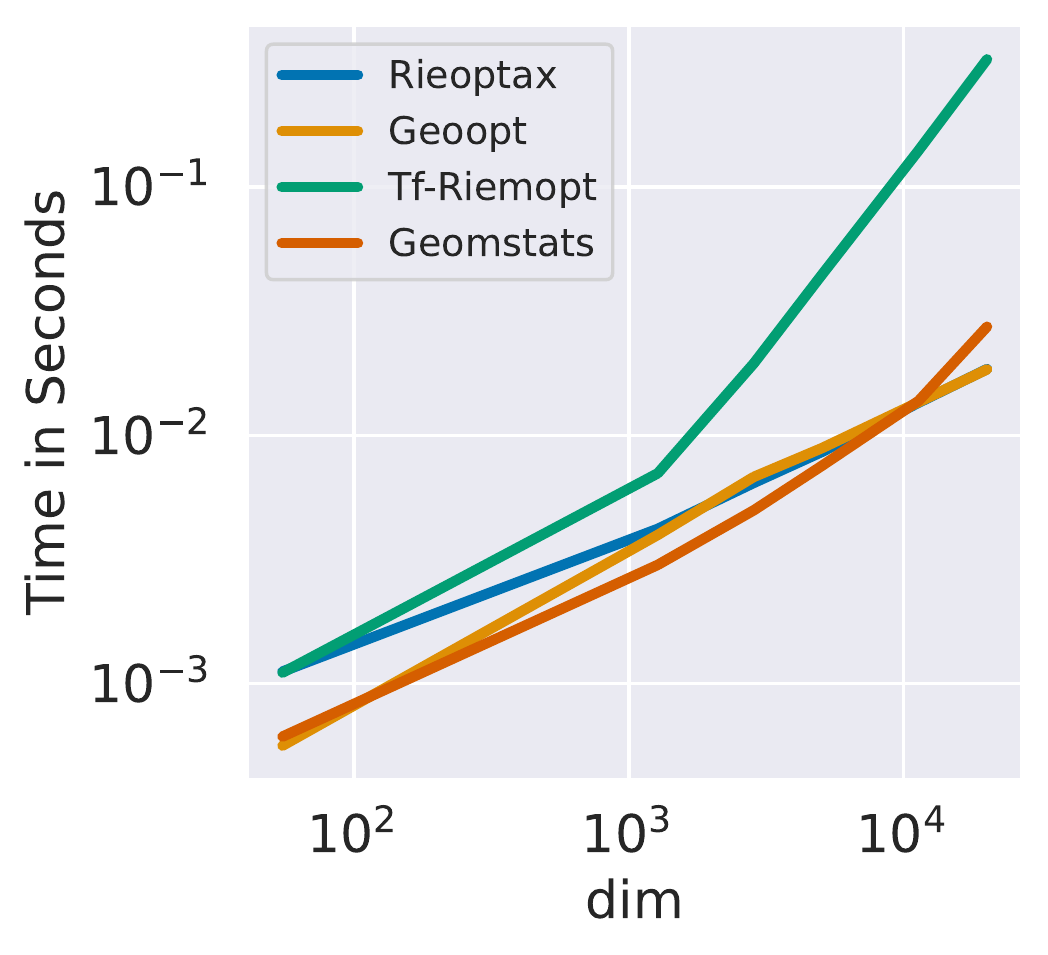}}

\caption{Benchmarking of Geometric Primitives on GPU.}\label{fig:GPUresults}
\end{figure}

\begin{itemize}
    \item  \textbf{Hypersphere:}  Hypersphere $\mathcal{S}(d)$ is supported in Geoopt, Tf-Riemopt, Geomstats, McTorch, and Pyamanopt. McTorch does not support the Exp and Log maps. On GPU, Geomstats raises an error. We benchmark for dimensions $ d \in \{50, 100,  500, 1000, 5000, 10000, 25000, 50000 \}$. 

    \item \textbf{Loretnz hyperboloid model:} The Loretnz hyperboloid model $\mathbb{H}(d)$ is supported in Geoopt, Tf-Riemopt, Geomstats, and Mctorch. While the Exp map is available in Mctorch, it does not implement the Log map. We benchmark for dimensions $d \in \{50, 100,  500, 1000, 5000,$ $10000, 25000, 50000\}$.

    \item \textbf{Grassmann:} Grassmann manifold $\mathcal{G}(m,r)$ is supported in Tf-Riemopt, Pymanopt, Geomstats. However, we notice that the logarithm map in Tf-Riemopt is incorrectly implemented 
    and Geomstats represents Grassmann elements in projector matrices form $\mathbf X \mathbf X^\top \in \mathbb R^{m \times m}$ instead of  $\mathbf X \in \mathbb R^{m \times r}$, which is prohibitively expensive. We thus exclude these two libraries from benchmarking. We benchmark for matrix sizes $(m,r) \in \{ (100,10) ,(500, 10), \allowbreak (750, 10), (1000, 10), (2000,10), (5000, 10) \}$.
    
    \item \textbf{SPD with affine-invariant metric:} SPD manifold ${\rm SPD}(m)$ with the affine-invariant metric is supported in Geoopt, Tf-Riemopt, Geomstats, and McTorch. McTorch, however, does not support the Exp and Log maps. We benchmark for matrix sizes
    $m \in \{10, 50, 75, 100,  150,  200 \}$.
\end{itemize}
Figures~\ref{fig:CPUresults} and~\ref{fig:GPUresults} present the timing results with CPU- and GPU-based computations, respectively. Overall, we observe that Rieoptax offers significant time improvements, especially on GPUs. For the SPDAffineInvariant case, Rieoptax is slightly slower than Geoopt because \texttt{eigh} which provides eigen decomposition is slightly slower in JAX compared to Pytorch. 
Given that JAX is a relatively new framework, we believe it would be faster even in this case in the near future.

\section{Example}

In this section, we show full code for the principal component analysis (PCA) by viewing it as optimization problem on Grassmann manifold. It can be found in Listing \ref{list:code}.

\begin{listing}
\begin{align*}
    \min_{\textbf{U} \in \mathcal{G}(m, r)} \frac{1}{n} \sum_{i=1}^{n} \norm{\textbf{z}_{i} - \textbf{U}\textbf{U}^T \textbf{z}_i}_{2}^2
\end{align*}
\begin{mdframed}
\caption{Principal component analysis (PCA) on the Grassmann manifold, where \texttt{init}, \texttt{Z} have to be provided. }
\label{list:code}
\begin{minted}{python}
from jax import jit, vmap
from jax.numpy.linalg import norm 
from rieoptax.core import rgrad, ManifoldArray
from rieoptax.geometry.grassmann import GrassmannCanonical
from rieoptax.optimizer.first_order import rsgd, dp_rsgd
from rieoptax.mechanism.gradient_perturbation import DP_RSGD_Mechanism
from rieoptax.optimizer.update import apply_updates

def fit(params, data, optimizer, epochs):
    @jit
    def step(params, opt_state, data):
    
        def cost(params, data):
            def _cost(params, data):
                return norm(data-params.value@(params.value.T @ data))**2
            return vmap(_cost, in_axes=(None,0))(params, data).mean()
            
        rgrads = rgrad(cost)(params, data)#calculate Riemannian grads
        updates, opt_state = optimizer.update(rgrads, opt_state, params)
        params = apply_updates(params, updates)#Uses Riemannian Exp
        return params, opt_state, loss_value
    
    opt_state = optimizer.init(params)
    for i in range(epochs):
        params, opt_state, loss_value = step(params, opt_state, data)

#initialization
U_init = ManifoldArray(value=init, manifold=GrassmannCanonical())
        
# non private PCA
lr, epochs = (0.1, 1000)
optimizer = rsgd(lr)
non_private_U = fit(U_init, Z, optimizer, epochs)

#(eps, delta) differentially private PCA
eps, delta, clip_norm = (1.0, 1e-6, 1.0)
sigma = DP_RGD_Mechanism(eps, delta, clip_norm, n)
optimizer = dp_rsgd(lr, sigma, clip_norm) 
private_U = fit(U_init, Z, optimizer, epochs)
\end{minted}
\end{mdframed}

\end{listing}

\section{Conclusion and future roadmap}
In this work, we present a Python library for (privacy-supported) Riemannian optimization, Rieoptax, and illustrate its efficacy on both CPU and GPU architectures. Our roadmap includes adding support for more manifold geometries, optimization algorithms, and a collection of example codes showcasing the usage of Rieoptax in various applications.


\bibliographystyle{plain}
\bibliography{references.bib}

\begin{thebibliography}{100}

\bibitem{abadi2016tensorflow}
M.~Abadi, A.~Agarwal, P.~Barham, E.~Brevdo, Z.~Chen, C.~Citro, G.~S. Corrado,
  A.~Davis, J.~Dean, M.~Devin, et~al.
\newblock Tensorflow: Large-scale machine learning on heterogeneous distributed
  systems.
\newblock In {\em USENIX Conference on Operating Systems Design and
  Implementation}, 2016.

\bibitem{abadi2016deep}
M.~Abadi, A.~Chu, I.~Goodfellow, H.~B. McMahan, I.~Mironov, K.~Talwar, and
  L.~Zhang.
\newblock Deep learning with differential privacy.
\newblock In {\em Proceedings of the 2016 ACM SIGSAC conference on computer and
  communications security}, pages 308--318, 2016.

\bibitem{abowd2018us}
J.~M. Abowd.
\newblock The {US Census Bureau} adopts differential privacy.
\newblock In {\em Proceedings of the 24th ACM SIGKDD International Conference
  on Knowledge Discovery \& Data Mining}, pages 2867--2867, 2018.

\bibitem{absil2007trust}
P.-A. Absil, C.~G. Baker, and K.~A. Gallivan.
\newblock Trust-region methods on {R}iemannian manifolds.
\newblock {\em Foundations of Computational Mathematics}, 7(3):303--330, 2007.

\bibitem{absil2009optimization}
P.-A. Absil, R.~Mahony, and R.~Sepulchre.
\newblock Optimization algorithms on matrix manifolds.
\newblock In {\em Optimization Algorithms on Matrix Manifolds}. Princeton
  University Press, 2009.

\bibitem{agarwal2021adaptive}
N.~Agarwal, N.~Boumal, B.~Bullins, and C.~Cartis.
\newblock Adaptive regularization with cubics on manifolds.
\newblock {\em Mathematical Programming}, 188(1):85--134, 2021.

\bibitem{ahn2020nesterov}
K.~Ahn and S.~Sra.
\newblock {From Nesterov’s estimate sequence to Riemannian acceleration}.
\newblock In {\em Conference on Learning Theory}, pages 84--118. PMLR, 2020.

\bibitem{alimisis2020continuous}
F.~Alimisis, A.~Orvieto, G.~B{\'e}cigneul, and A.~Lucchi.
\newblock A continuous-time perspective for modeling acceleration in
  {R}iemannian optimization.
\newblock In {\em International Conference on Artificial Intelligence and
  Statistics}, pages 1297--1307. PMLR, 2020.

\bibitem{altschuler2021averaging}
J.~Altschuler, S.~Chewi, P.~R. Gerber, and A.~Stromme.
\newblock Averaging on the {Bures-Wasserstein} manifold: dimension-free
  convergence of gradient descent.
\newblock {\em Advances in Neural Information Processing Systems},
  34:22132--22145, 2021.

\bibitem{apple2017learning}
D.~Apple.
\newblock Learning with privacy at scale.
\newblock {\em Apple Machine Learning Journal}, 1(8), 2017.

\bibitem{arsigny2007geometric}
V.~Arsigny, P.~Fillard, X.~Pennec, and N.~Ayache.
\newblock Geometric means in a novel vector space structure on symmetric
  positive-definite matrices.
\newblock {\em SIAM journal on matrix analysis and applications},
  29(1):328--347, 2007.

\bibitem{deepmind2020jax}
I.~Babuschkin, K.~Baumli, A.~Bell, S.~Bhupatiraju, J.~Bruce, P.~Buchlovsky,
  D.~Budden, T.~Cai, A.~Clark, I.~Danihelka, C.~Fantacci, J.~Godwin, C.~Jones,
  R.~Hemsley, T.~Hennigan, M.~Hessel, S.~Hou, S.~Kapturowski, T.~Keck,
  I.~Kemaev, M.~King, M.~Kunesch, L.~Martens, H.~Merzic, V.~Mikulik, T.~Norman,
  J.~Quan, G.~Papamakarios, R.~Ring, F.~Ruiz, A.~Sanchez, R.~Schneider,
  E.~Sezener, S.~Spencer, S.~Srinivasan, L.~Wang, W.~Stokowiec, and F.~Viola.
\newblock The {D}eep{M}ind {JAX} {E}cosystem, 2020.

\bibitem{balle2022reconstructing}
B.~Balle, G.~Cherubin, and J.~Hayes.
\newblock Reconstructing training data with informed adversaries.
\newblock {\em arXiv preprint arXiv:2201.04845}, 2022.

\bibitem{baydin2018automatic}
A.~G. Baydin, B.~A. Pearlmutter, A.~A. Radul, and J.~M. Siskind.
\newblock Automatic differentiation in machine learning: a survey.
\newblock {\em Journal of Marchine Learning Research}, 18:1--43, 2018.

\bibitem{becigneul2018riemannian}
G.~Becigneul and O.-E. Ganea.
\newblock Riemannian adaptive optimization methods.
\newblock In {\em International Conference on Learning Representations}, 2019.

\bibitem{bendokat2020grassmann}
T.~Bendokat, R.~Zimmermann, and P.-A. Absil.
\newblock A {G}rassmann manifold handbook: Basic geometry and computational
  aspects.
\newblock {\em arXiv preprint arXiv:2011.13699}, 2020.

\bibitem{bergmann2022manopt}
R.~Bergmann.
\newblock Manopt. jl: Optimization on manifolds in julia.
\newblock {\em Journal of Open Source Software}, 7(70):3866, 2022.

\bibitem{bezanson2017julia}
J.~Bezanson, A.~Edelman, S.~Karpinski, and V.~B. Shah.
\newblock Julia: A fresh approach to numerical computing.
\newblock {\em SIAM review}, 59(1):65--98, 2017.

\bibitem{bhatia2009positive}
R.~Bhatia.
\newblock Positive definite matrices.
\newblock In {\em Positive Definite Matrices}. Princeton university press,
  2009.

\bibitem{bonnabel2013stochastic}
S.~Bonnabel.
\newblock Stochastic gradient descent on {R}iemannian manifolds.
\newblock {\em IEEE Transactions on Automatic Control}, 58(9):2217--2229, 2013.

\bibitem{boumal2022intromanifolds}
N.~Boumal.
\newblock An introduction to optimization on smooth manifolds.
\newblock To appear with Cambridge University Press, Jun 2022.

\bibitem{boumal2011rtrmc}
N.~Boumal and P.-a. Absil.
\newblock {RTRMC: A Riemannian trust-region method for low-rank matrix
  completion}.
\newblock {\em Advances in neural information processing systems}, 24, 2011.

\bibitem{boumal2014manopt}
N.~Boumal, B.~Mishra, P.-A. Absil, and R.~Sepulchre.
\newblock Manopt, a {M}atlab toolbox for optimization on manifolds.
\newblock {\em The Journal of Machine Learning Research}, 15(1):1455--1459,
  2014.

\bibitem{jax2018github}
J.~Bradbury, R.~Frostig, P.~Hawkins, M.~J. Johnson, C.~Leary, D.~Maclaurin,
  G.~Necula, A.~Paszke, J.~Vander{P}las, S.~Wanderman-{M}ilne, and Q.~Zhang.
\newblock {JAX}: composable transformations of {P}ython+{N}um{P}y programs,
  2018.

\bibitem{bu2021fast}
Z.~Bu, S.~Gopi, J.~Kulkarni, Y.~T. Lee, H.~Shen, and U.~Tantipongpipat.
\newblock Fast and memory efficient differentially private-sgd via jl
  projections.
\newblock {\em Advances in Neural Information Processing Systems},
  34:19680--19691, 2021.

\bibitem{10.5555/3361338.3361358}
N.~Carlini, C.~Liu, U.~Erlingsson, J.~Kos, and D.~Song.
\newblock The secret sharer: Evaluating and testing unintended memorization in
  neural networks.
\newblock In {\em Proceedings of the 28th USENIX Conference on Security
  Symposium}, SEC'19, page 267–284, USA, 2019. USENIX Association.

\bibitem{chakraborty2019statistics}
R.~Chakraborty and B.~C. Vemuri.
\newblock Statistics on the {S}tiefel manifold: theory and applications.
\newblock {\em The Annals of Statistics}, 47(1):415--438, 2019.

\bibitem{chaudhuri2011differentially}
K.~Chaudhuri, C.~Monteleoni, and A.~D. Sarwate.
\newblock Differentially private empirical risk minimization.
\newblock {\em Journal of Machine Learning Research}, 12(3), 2011.

\bibitem{chewi2020gradient}
S.~Chewi, T.~Maunu, P.~Rigollet, and A.~J. Stromme.
\newblock Gradient descent algorithms for {Bures-Wasserstein} barycenters.
\newblock In {\em Conference on Learning Theory}, pages 1276--1304. PMLR, 2020.

\bibitem{cuturi2022optimal}
M.~Cuturi, L.~Meng-Papaxanthos, Y.~Tian, C.~Bunne, G.~Davis, and O.~Teboul.
\newblock Optimal transport tools ({OTT}): A jax toolbox for all things
  {W}asserstein.
\newblock {\em arXiv preprint arXiv:2201.12324}, 2022.

\bibitem{ding2017collecting}
B.~Ding, J.~Kulkarni, and S.~Yekhanin.
\newblock Collecting telemetry data privately.
\newblock {\em Advances in Neural Information Processing Systems}, 30, 2017.

\bibitem{dwork2008differential}
C.~Dwork.
\newblock Differential privacy: A survey of results.
\newblock In {\em International conference on theory and applications of models
  of computation}, pages 1--19. Springer, 2008.

\bibitem{dwork2006our}
C.~Dwork, K.~Kenthapadi, F.~McSherry, I.~Mironov, and M.~Naor.
\newblock Our data, ourselves: Privacy via distributed noise generation.
\newblock In {\em Annual international conference on the theory and
  applications of cryptographic techniques}, pages 486--503. Springer, 2006.

\bibitem{dwork2006calibrating}
C.~Dwork, F.~McSherry, K.~Nissim, and A.~Smith.
\newblock Calibrating noise to sensitivity in private data analysis.
\newblock In {\em Theory of cryptography conference}, pages 265--284. Springer,
  2006.

\bibitem{dwork2014algorithmic}
C.~Dwork, A.~Roth, et~al.
\newblock The algorithmic foundations of differential privacy.
\newblock {\em Foundations and Trends{\textregistered} in Theoretical Computer
  Science}, 9(3--4):211--407, 2014.

\bibitem{edelman1998geometry}
A.~Edelman, T.~A. Arias, and S.~T. Smith.
\newblock The geometry of algorithms with orthogonality constraints.
\newblock {\em SIAM journal on Matrix Analysis and Applications},
  20(2):303--353, 1998.

\bibitem{erlingsson2014rappor}
{\'U}.~Erlingsson, V.~Pihur, and A.~Korolova.
\newblock Rappor: {R}andomized aggregatable privacy-preserving ordinal
  response.
\newblock In {\em Proceedings of the 2014 ACM SIGSAC conference on computer and
  communications security}, pages 1054--1067, 2014.

\bibitem{finzi2021practical}
M.~Finzi, M.~Welling, and A.~G. Wilson.
\newblock A practical method for constructing equivariant multilayer
  perceptrons for arbitrary matrix groups.
\newblock In {\em International Conference on Machine Learning}, pages
  3318--3328. PMLR, 2021.

\bibitem{frechet1948elements}
M.~Fr{\'e}chet.
\newblock Les {\'e}l{\'e}ments al{\'e}atoires de nature quelconque dans un
  espace distanci{\'e}.
\newblock In {\em Annales de l'institut Henri Poincar{\'e}}, volume~10, pages
  215--310, 1948.

\bibitem{brax2021github}
C.~D. Freeman, E.~Frey, A.~Raichuk, S.~Girgin, I.~Mordatch, and O.~Bachem.
\newblock Brax - a differentiable physics engine for large scale rigid body
  simulation, 2021.

\bibitem{frostig2018compiling}
R.~Frostig, M.~J. Johnson, and C.~Leary.
\newblock Compiling machine learning programs via high-level tracing.
\newblock {\em Systems for Machine Learning}, 4(9), 2018.

\bibitem{gallier2020differential}
J.~Gallier and J.~Quaintance.
\newblock {\em Differential geometry and {L}ie groups: a computational
  perspective}, volume~12.
\newblock Springer Nature, 2020.

\bibitem{ganea2018hyperbolic}
O.~Ganea, G.~B{\'e}cigneul, and T.~Hofmann.
\newblock Hyperbolic neural networks.
\newblock {\em Advances in neural information processing systems}, 31, 2018.

\bibitem{jraph2020github}
J.~Godwin*, T.~Keck*, P.~Battaglia, V.~Bapst, T.~Kipf, Y.~Li, K.~Stachenfeld,
  P.~Veli\v{c}kovi\'{c}, and A.~Sanchez-Gonzalez.
\newblock {J}raph: {A} library for graph neural networks in jax., 2020.

\bibitem{goodfellow2015efficient}
I.~Goodfellow.
\newblock Efficient per-example gradient computations.
\newblock {\em arXiv preprint arXiv:1510.01799}, 2015.

\bibitem{50530}
Google.
\newblock Xla : Compiling machine learning for peak performance, 2020.

\bibitem{han2021improved}
A.~Han and J.~Gao.
\newblock {Improved variance reduction methods for Riemannian non-convex
  optimization}.
\newblock {\em IEEE Transactions on Pattern Analysis and Machine Intelligence},
  2021.

\bibitem{hanmomentum2021}
A.~Han and J.~Gao.
\newblock Riemannian stochastic recursive momentum method for non-convex
  optimization.
\newblock In {\em International Joint Conference on Artificial Intelligence},
  pages 2505--2511, 8 2021.

\bibitem{han2022differentially}
A.~Han, B.~Mishra, P.~Jawanpuria, and J.~Gao.
\newblock Differentially private {R}iemannian optimization.
\newblock {\em arXiv preprint arXiv:2205.09494}, 2022.

\bibitem{han2022riemannian}
A.~Han, B.~Mishra, P.~Jawanpuria, and J.~Gao.
\newblock Riemannian accelerated gradient methods via extrapolation.
\newblock {\em arXiv preprint arXiv:2208.06619}, 2022.

\bibitem{han2022riemannianblockSPD}
A.~Han, B.~Mishra, P.~Jawanpuria, and J.~Gao.
\newblock Riemannian block {SPD} coupling manifold and its application to
  optimal transport.
\newblock {\em arXiv preprint arXiv:2201.12933}, 2022.

\bibitem{han2021riemannian}
A.~Han, B.~Mishra, P.~K. Jawanpuria, and J.~Gao.
\newblock On {Riemannian} optimization over positive definite matrices with the
  {Bures-Wasserstein} geometry.
\newblock {\em Advances in Neural Information Processing Systems},
  34:8940--8953, 2021.

\bibitem{harris2020array}
C.~R. Harris, K.~J. Millman, S.~J. Van Der~Walt, R.~Gommers, P.~Virtanen,
  D.~Cournapeau, E.~Wieser, J.~Taylor, S.~Berg, N.~J. Smith, et~al.
\newblock {Array programming with numpy}.
\newblock {\em Nature}, 585(7825):357--362, 2020.

\bibitem{flax2020github}
J.~Heek, A.~Levskaya, A.~Oliver, M.~Ritter, B.~Rondepierre, A.~Steiner, and
  M.~van {Z}ee.
\newblock {F}lax: A neural network library and ecosystem for {JAX}, 2020.

\bibitem{helgason1979differential}
S.~Helgason.
\newblock {\em Differential geometry, Lie groups, and symmetric spaces}.
\newblock Academic press, 1979.

\bibitem{haiku2020github}
T.~Hennigan, T.~Cai, T.~Norman, and I.~Babuschkin.
\newblock {H}aiku: {S}onnet for {JAX}, 2020.

\bibitem{hosseini2020alternative}
R.~Hosseini and S.~Sra.
\newblock {An alternative to EM for Gaussian mixture models: batch and
  stochastic Riemannian optimization}.
\newblock {\em Mathematical programming}, 181(1):187--223, 2020.

\bibitem{huang2015broyden}
W.~Huang, K.~A. Gallivan, and P.-A. Absil.
\newblock {A Broyden class of quasi-Newton methods for Riemannian
  optimization}.
\newblock {\em SIAM Journal on Optimization}, 25(3):1660--1685, 2015.

\bibitem{huang2016riemannian}
W.~Huang, K.~A. Gallivan, A.~Srivastava, and P.-A. Absil.
\newblock Riemannian optimization for registration of curves in elastic shape
  analysis.
\newblock {\em Journal of Mathematical Imaging and Vision}, 54(3):320--343,
  2016.

\bibitem{huang2017riemannian}
Z.~Huang and L.~Van~Gool.
\newblock A {R}iemannian network for {SPD} matrix learning.
\newblock In {\em Thirty-first AAAI conference on artificial intelligence},
  2017.

\bibitem{huang2017deep}
Z.~Huang, C.~Wan, T.~Probst, and L.~Van~Gool.
\newblock Deep learning on {Lie} groups for skeleton-based action recognition.
\newblock In {\em Proceedings of the IEEE conference on computer vision and
  pattern recognition}, pages 6099--6108, 2017.

\bibitem{huang2018building}
Z.~Huang, J.~Wu, and L.~Van~Gool.
\newblock Building deep networks on {Grassmann} manifolds.
\newblock In {\em Proceedings of the AAAI Conference on Artificial
  Intelligence}, 2018.

\bibitem{kachan2020persistent}
O.~Kachan.
\newblock Persistent homology-based projection pursuit.
\newblock In {\em Proceedings of the IEEE/CVF Conference on Computer Vision and
  Pattern Recognition Workshops}, pages 856--857, 2020.

\bibitem{kasai2019riemannian}
H.~Kasai, P.~Jawanpuria, and B.~Mishra.
\newblock Riemannian adaptive stochastic gradient algorithms on matrix
  manifolds.
\newblock In {\em International Conference on Machine Learning}, pages
  3262--3271. PMLR, 2019.

\bibitem{kasai2018riemannian}
H.~Kasai, H.~Sato, and B.~Mishra.
\newblock Riemannian stochastic recursive gradient algorithm.
\newblock In {\em International Conference on Machine Learning}, pages
  2516--2524. PMLR, 2018.

\bibitem{kendall1984shape}
D.~G. Kendall.
\newblock Shape manifolds, procrustean metrics, and complex projective spaces.
\newblock {\em Bulletin of the London mathematical society}, 16(2):81--121,
  1984.

\bibitem{kendall1989survey}
D.~G. Kendall.
\newblock A survey of the statistical theory of shape.
\newblock {\em Statistical Science}, 4(2):87--99, 1989.

\bibitem{kidger2021on}
P.~Kidger.
\newblock {\em {O}n {N}eural {D}ifferential {E}quations}.
\newblock PhD thesis, University of Oxford, 2021.

\bibitem{kidger2021equinox}
P.~Kidger and C.~Garcia.
\newblock {E}quinox: neural networks in {JAX} via callable {P}y{T}rees and
  filtered transformations.
\newblock {\em Differentiable Programming workshop at Neural Information
  Processing Systems 2021}, 2021.

\bibitem{kochurov2020geoopt}
M.~Kochurov, R.~Karimov, and S.~Kozlukov.
\newblock {Geoopt: Riemannian optimization in PyTorch}.
\newblock {\em arXiv preprint arXiv:2005.02819}, 2020.

\bibitem{blackjax2020github}
J.~Lao and R.~Louf.
\newblock {B}lackjax: A sampling library for {JAX}, 2020.

\bibitem{le2015tiny}
Y.~Le and X.~Yang.
\newblock Tiny imagenet visual recognition challenge.
\newblock {\em CS 231N}, 7(7):3, 2015.

\bibitem{lee2021scaling}
J.~Lee and D.~Kifer.
\newblock Scaling up differentially private deep learning with fast per-example
  gradient clipping.
\newblock {\em Proceedings on Privacy Enhancing Technologies}, 2021(1), 2021.

\bibitem{lee2006riemannian}
J.~M. Lee.
\newblock {\em Riemannian manifolds: an introduction to curvature}, volume 176.
\newblock Springer Science \& Business Media, 2006.

\bibitem{li2022stochastic}
J.~Li, K.~Balasubramanian, and S.~Ma.
\newblock Stochastic zeroth-order {R}iemannian derivative estimation and
  optimization.
\newblock {\em Mathematics of Operations Research}, 2022.

\bibitem{liu2017accelerated}
Y.~Liu, F.~Shang, J.~Cheng, H.~Cheng, and L.~Jiao.
\newblock Accelerated first-order methods for geodesically convex optimization
  on {R}iemannian manifolds.
\newblock {\em Advances in Neural Information Processing Systems}, 30, 2017.

\bibitem{maclaurin2015autograd}
D.~Maclaurin, D.~Duvenaud, and R.~P. Adams.
\newblock Autograd: Effortless gradients in numpy.
\newblock In {\em ICML 2015 AutoML workshop}, 2015.

\bibitem{meghwanshi2018mctorch}
M.~Meghwanshi, P.~Jawanpuria, A.~Kunchukuttan, H.~Kasai, and B.~Mishra.
\newblock {McTorch, a manifold optimization library for deep learning}.
\newblock {\em arXiv preprint arXiv:1810.01811}, 2018.

\bibitem{miolane2020geomstats}
N.~Miolane, N.~Guigui, A.~Le~Brigant, J.~Mathe, B.~Hou, Y.~Thanwerdas,
  S.~Heyder, O.~Peltre, N.~Koep, H.~Zaatiti, et~al.
\newblock {Geomstats: a Python package for Riemannian geometry in machine
  learning}.
\newblock {\em Journal of Machine Learning Research}, 21(223):1--9, 2020.

\bibitem{miolane2017template}
N.~Miolane, S.~Holmes, and X.~Pennec.
\newblock Template shape estimation: correcting an asymptotic bias.
\newblock {\em SIAM Journal on Imaging Sciences}, 10(2):808--844, 2017.

\bibitem{mironov2017renyi}
I.~Mironov.
\newblock R{\'e}nyi differential privacy.
\newblock In {\em 2017 IEEE 30th computer security foundations symposium
  (CSF)}, pages 263--275. IEEE, 2017.

\bibitem{mishra2021manifold}
B.~Mishra, N.~Satyadev, H.~Kasai, and P.~Jawanpuria.
\newblock Manifold optimization for non-linear optimal transport problems.
\newblock {\em arXiv preprint arXiv:2103.00902}, 2021.

\bibitem{near2018differential}
J.~Near.
\newblock Differential privacy at scale: {Uber and Berkeley} collaboration.
\newblock In {\em Enigma 2018 (Enigma 2018)}, 2018.

\bibitem{nguyen2019neural}
X.~S. Nguyen, L.~Brun, O.~L{\'e}zoray, and S.~Bougleux.
\newblock A neural network based on spd manifold learning for skeleton-based
  hand gesture recognition.
\newblock In {\em Proceedings of the IEEE/CVF Conference on Computer Vision and
  Pattern Recognition}, pages 12036--12045, 2019.

\bibitem{nickel2017poincare}
M.~Nickel and D.~Kiela.
\newblock Poincar{\'e} embeddings for learning hierarchical representations.
\newblock {\em Advances in neural information processing systems}, 30, 2017.

\bibitem{nickel2018learning}
M.~Nickel and D.~Kiela.
\newblock Learning continuous hierarchies in the {L}orentz model of hyperbolic
  geometry.
\newblock In {\em International Conference on Machine Learning}, pages
  3779--3788. PMLR, 2018.

\bibitem{paszke2019pytorch}
A.~Paszke, S.~Gross, F.~Massa, A.~Lerer, J.~Bradbury, G.~Chanan, T.~Killeen,
  Z.~Lin, N.~Gimelshein, L.~Antiga, et~al.
\newblock Pytorch: An imperative style, high-performance deep learning library.
\newblock {\em Advances in neural information processing systems}, 32, 2019.

\bibitem{pennec2006riemannian}
X.~Pennec, P.~Fillard, and N.~Ayache.
\newblock A {R}iemannian framework for tensor computing.
\newblock {\em International Journal of computer vision}, 66(1):41--66, 2006.

\bibitem{qi2010riemannian}
C.~Qi, K.~A. Gallivan, and P.-A. Absil.
\newblock Riemannian {BFGS} algorithm with applications.
\newblock In {\em Recent Advances in Optimization and its Applications in
  Engineering}, pages 183--192. Springer, 2010.

\bibitem{qi2021transductive}
G.~Qi, H.~Yu, Z.~Lu, and S.~Li.
\newblock Transductive few-shot classification on the oblique manifold.
\newblock In {\em Proceedings of the IEEE/CVF International Conference on
  Computer Vision}, pages 8412--8422, 2021.

\bibitem{rahman2018membership}
M.~A. Rahman, T.~Rahman, R.~Lagani{\`e}re, N.~Mohammed, and Y.~Wang.
\newblock Membership inference attack against differentially private deep
  learning model.
\newblock {\em Trans. Data Priv.}, 11(1):61--79, 2018.

\bibitem{reimherr2021differential}
M.~Reimherr, K.~Bharath, and C.~Soto.
\newblock Differential privacy over {R}iemannian manifolds.
\newblock {\em Advances in Neural Information Processing Systems},
  34:12292--12303, 2021.

\bibitem{fedjax2021}
J.~H. Ro, A.~T. Suresh, and K.~Wu.
\newblock {F}ed{JAX}: Federated learning simulation with {JAX}.
\newblock {\em arXiv preprint arXiv:2108.02117}, 2021.

\bibitem{rochette2019efficient}
G.~Rochette, A.~Manoel, and E.~W. Tramel.
\newblock Efficient per-example gradient computations in convolutional neural
  networks.
\newblock {\em arXiv preprint arXiv:1912.06015}, 2019.

\bibitem{sablayrolles2019white}
A.~Sablayrolles, M.~Douze, C.~Schmid, Y.~Ollivier, and H.~J{\'e}gou.
\newblock White-box vs black-box: Bayes optimal strategies for membership
  inference.
\newblock In {\em International Conference on Machine Learning}, pages
  5558--5567. PMLR, 2019.

\bibitem{sato2019riemannian}
H.~Sato, H.~Kasai, and B.~Mishra.
\newblock Riemannian stochastic variance reduced gradient algorithm with
  retraction and vector transport.
\newblock {\em SIAM Journal on Optimization}, 29(2):1444--1472, 2019.

\bibitem{jaxmd2020}
S.~S. Schoenholz and E.~D. Cubuk.
\newblock Jax m.d. a framework for differentiable physics.
\newblock In {\em Advances in Neural Information Processing Systems},
  volume~33. Curran Associates, Inc., 2020.

\bibitem{selig2005geometric}
J.~M. Selig.
\newblock {\em Geometric fundamentals of robotics}, volume 128.
\newblock Springer, 2005.

\bibitem{shi2021coupling}
D.~Shi, J.~Gao, X.~Hong, S.~Boris~Choy, and Z.~Wang.
\newblock Coupling matrix manifolds assisted optimization for optimal transport
  problems.
\newblock {\em Machine Learning}, 110(3):533--558, 2021.

\bibitem{smirnov2021tensorflow}
O.~Smirnov.
\newblock {TensorFlow RiemOpt: a library for optimization on Riemannian
  manifolds}.
\newblock {\em arXiv preprint arXiv:2105.13921}, 2021.

\bibitem{sola2018micro}
J.~Sola, J.~Deray, and D.~Atchuthan.
\newblock A micro {L}ie theory for state estimation in robotics.
\newblock {\em arXiv preprint arXiv:1812.01537}, 2018.

\bibitem{song2013stochastic}
S.~Song, K.~Chaudhuri, and A.~D. Sarwate.
\newblock Stochastic gradient descent with differentially private updates.
\newblock In {\em 2013 IEEE global conference on signal and information
  processing}, pages 245--248. IEEE, 2013.

\bibitem{srivastava2010shape}
A.~Srivastava, E.~Klassen, S.~H. Joshi, and I.~H. Jermyn.
\newblock Shape analysis of elastic curves in {E}uclidean spaces.
\newblock {\em IEEE transactions on pattern analysis and machine intelligence},
  33(7):1415--1428, 2010.

\bibitem{subramani2021enabling}
P.~Subramani, N.~Vadivelu, and G.~Kamath.
\newblock Enabling fast differentially private {SGD} via just-in-time
  compilation and vectorization.
\newblock {\em Advances in Neural Information Processing Systems},
  34:26409--26421, 2021.

\bibitem{thanwerdas2021n}
Y.~Thanwerdas and X.~Pennec.
\newblock O(n)-invariant {R}iemannian metrics on {SPD} matrices.
\newblock {\em arXiv preprint arXiv:2109.05768}, 2021.

\bibitem{JMLR:v17:16-177}
J.~Townsend, N.~Koep, and S.~Weichwald.
\newblock Pymanopt: A {Python} toolbox for optimization on manifolds using
  automatic differentiation.
\newblock {\em Journal of Machine Learning Research}, 17(137):1--5, 2016.

\bibitem{ungar2008analytic}
A.~A. Ungar.
\newblock {\em Analytic hyperbolic geometry and {A}lbert {E}instein’s special
  theory of relativity}.
\newblock World Scientific, 2008.

\bibitem{ungar2008gyrovector}
A.~A. Ungar.
\newblock A gyrovector space approach to hyperbolic geometry.
\newblock {\em Synthesis Lectures on Mathematics and Statistics}, 1(1):1--194,
  2008.

\bibitem{utpala2022differentially}
S.~Utpala, P.~Vepakomma, and N.~Miolane.
\newblock Differentially private {F}r\'echet mean on the manifold of symmetric
  positive definite ({SPD}) matrices.
\newblock {\em arXiv preprint arXiv:2208.04245}, 2022.

\bibitem{wang2019subsampled}
Y.-X. Wang, B.~Balle, and S.~P. Kasiviswanathan.
\newblock Subsampled r{\'e}nyi differential privacy and analytical moments
  accountant.
\newblock In {\em The 22nd International Conference on Artificial Intelligence
  and Statistics}, pages 1226--1235. PMLR, 2019.

\bibitem{zhang2016riemannian}
H.~Zhang, S.~J~Reddi, and S.~Sra.
\newblock {Riemannian SVRG: Fast stochastic optimization on Riemannian
  manifolds}.
\newblock {\em Advances in Neural Information Processing Systems}, 29, 2016.

\bibitem{zhang2016first}
H.~Zhang and S.~Sra.
\newblock First-order methods for geodesically convex optimization.
\newblock In {\em Conference on Learning Theory}, pages 1617--1638. PMLR, 2016.

\bibitem{zhang2018estimate}
H.~Zhang and S.~Sra.
\newblock An estimate sequence for geodesically convex optimization.
\newblock In {\em Conference On Learning Theory}, pages 1703--1723. PMLR, 2018.

\bibitem{zhou2019faster}
P.~Zhou, X.-T. Yuan, and J.~Feng.
\newblock {Faster first-order methods for stochastic non-convex optimization on
  Riemannian manifolds}.
\newblock In {\em The 22nd International Conference on Artificial Intelligence
  and Statistics}, pages 138--147. PMLR, 2019.

\bibitem{zhu2019deep}
L.~Zhu, Z.~Liu, and S.~Han.
\newblock Deep leakage from gradients.
\newblock {\em Advances in neural information processing systems}, 32, 2019.

\end{thebibliography}

\end{document}